\documentclass[a4paper]{article}
\usepackage{amsmath,amssymb,amsthm,amscd}

\newtheorem{Thm}{Theorem}[subsection]
\newtheorem{Lem}[Thm]{Lemma}
\newtheorem{Cor}[Thm]{Corollary}
\newtheorem{Prop}[Thm]{Proposition}
\newtheorem{Def}[Thm]{Definition}
\newtheorem{Rem}[Thm]{Remark}
\newtheorem{Exam}[Thm]{Example}

\makeatletter
\@addtoreset{equation}{section}
\makeatother

\newcommand{\R}{\mathbb{R}}
\newcommand{\He}{\mathbb{H}}
\newcommand{\Ro}{\overline{\R}}
\newcommand{\N}{\mathbb{N}}
\newcommand{\M}{\mathbb{M}}
\newcommand{\Z}{\mathbb{Z}}
\newcommand{\C}{\mathbb{C}}

\newcommand{\K}{\mathcal{K}}
\newcommand{\An}{\mathcal{A}}
\newcommand{\Sch}{\mathcal{S}}
\newcommand{\E}{\mathcal{E}}
\newcommand{\Li}{\mathcal{L}}
\newcommand{\Hi}{\mathcal{H}}
\newcommand{\Fu}{\mathcal{F}}
\newcommand{\We}{\mathcal{W}}
\newcommand{\Mu}{\mathcal{M}}
\newcommand{\oneh}{\frac{1}{2}}
\newcommand{\Rb}{\overline{\R}}

\newcommand{\dbar}{{d\hspace{-0,05cm}\bar{}\hspace{0,05cm}}}
\newcommand{\op}{\textup{Op}}

\newcommand{\sop}{\textup{op}}

\newcommand{\clw}{\textup{cl}}
\newcommand{\loc}{\textup{loc}}
\newcommand{\cone}{\textup{cone}}

\newcommand{\la}{\langle}
\newcommand{\ra}{\rangle}
\newcommand{\lr}{\textup{(}}
\newcommand{\rr}{\textup{)}}

\newcommand{\imb}{\textup{Im}}
\newcommand{\ind}{\textup{ind}}

\newcommand{\reb}{\textup{Re}}

\newcommand{\vp}{\varphi}

\newcommand{\td}{\tilde}

\newcommand{\Diff}{\textup{Diff}}

\begin{document}

\title{Operators with Corner-degenerate Symbols}
\date{}
\author{Jamil Abed \and B.-W. Schulze}

\pagestyle{headings}

\maketitle

\begin{abstract}
We establish elements of a new approch to ellipticity and parametrices within operator algebras on a manifold with higher singularities, only based on some general axiomatic requirements on parameter-dependent operators in suitable scales of spaces. The idea is to model an iterative process with new generations of parameter-dependent operator theories, together with new scales of spaces that satisfy analogous requirements as the original ones, now on a corresponding higher level.\\ 
The ``full'' calculus is voluminous; so we content ourselves here with some typical aspects such as symbols in terms of order reducing families, classes of relevant examples, and operators near a corner point.
\end{abstract}
%\setcounter{footnote}{ }
%\footnotetext{2000 \textit{Mathematics Subject Classification}. Primary 35S35; Secondary 35J70.}
%\footnotetext{\textit{Key words and phrases}. pseudo-differential operators, weighted spaces, edge- and corner-degenerate symbols, ellipticity near conical exits to infinity}

\tableofcontents

\section*{Introduction}
\addcontentsline{toc}{section}{Introduction}
\markboth{INTRODUCTION}{INTRODUCTION}
This paper is aimed at studying operators with certain degenerate operator-valued amplitude functions, motivated by the iterative calculus of pseudo-differential operators on manifolds with higher singularities. Here, in contrast to  \cite{Schu27}, \cite{Schu49}, we develop the aspect of symbols, based  on ``abstract'' reductions of orders which makes the approch transparent from a new point of view. To illustrate the idea, let us first consider, for example, the Laplacian on a manifold with conical singularities (say, without boundary). In this case the ellipticity does not only refer to the ``standard'' principal homogeneous symbol but 
also to the so-called conormal symbol. The latter one, contributed by the conical point, is 
operator-valued and singles out the weights in Sobolev spaces, where the operator has the Fredholm property. 
Another example of ellipticity with different principal symbolic components is the case of boundary
value problems. The boundary (say, smooth), interpreted as an edge, contributes the 
operator-valued boundary (or edge) symbol which is responsible for the nature of boundary conditions
(for instance, of Dirichlet or Neumann type in the case of the Laplacian). In general, if the 
configuration has polyhedral singularities of order $k$, we have to expect a principal symbolic 
hierarchy of length $k+1$, with components contributed by the various strata.
In order to characterise the solvability of elliptic equations, especially, the regularity of 
solutions in suitable scales of spaces, it is adequate to embed the problem in a pseudo-differential 
calculus, and to construct a parametrix. For higher singularities this is a program of tremendous 
complexity. It is therefore advisable to organise general elements of the calculus by means of an axiomatic framework which contains the typical features, such as the cone- or edge-degenerate 
behaviour of symbols but ignores the (in general) huge tail of $k-1$ iterative steps to reach 
the singularity level $k$.\\
The ``concrete'' (pseudo-differential) calculus of operators on manifolds with conical or 
edge  singularities may be found in several papers and monographs, see, for instance, 
\cite{Remp7}, \cite{Schu2}, \cite{Schu32}, \cite{Egor1}. Operators on manifolds of singularity 
order 2 are studied in \cite{Schu29}, \cite{Schu27}, \cite{Mani2}, \cite{Haru11}. 
Theories of that kind are also possible for boundary value problems with the 
transmission property at the (smooth part of the) boundary, see, for instance, \cite{Remp11}, 
\cite{Kapa10}, \cite{Haru13}. This is useful in numerous applications, for instance, 
to models of elasticity or crack theory, see \cite{Kapa10}, \cite{Haru12}. Elements of operator structures on manifolds with higher singularities are developed, for instance, in \cite{Schu25}, \cite{Calv2}. The nature of such theories depends very much on specific assumptions on the degeneracy of the involved symbols. There are worldwide different schools studying operators on singular manifolds, partly motivated by problems of geometry, index theory, and topology, see, for instance, Melrose \cite{Melr2}, Melrose and Piazza \cite{Melr5}, Nistor \cite{Nist3}, Nazaikinskij, Savin, Sternin \cite{Naza18}, \cite{Naza19}, \cite{Naza17}, and many others. We do not study here operators of ``multi-Fuchs" type, often associated with ``corner manifolds". Our operators are of a rather different behaviour with respect to the degeneracy of symbols. Nevertheless the various theories have intersections and common sources, see the paper of Kondratiev \cite{Kond1} or papers and monographs of other representatives of a corresponding Russian school, see, for instance, 
\cite{Plam1}, \cite{Plam2}.\\
Let us briefly recall a few basic facts on operators on manifolds with conical singularities 
or edges.\\
Let $M$ be a manifold with conical singularity $v\in M$, i.e., $M \setminus \{ v \}$ is smooth, 
and $M$ is close to $v$ modelled on a cone $X^\Delta:=(\Rb_+\times X)/(\{0\}\times X)$ with base 
$X$, where $X$ is a closed compact $C^\infty$ manifold. We then have differential operators 
of order $\mu\in\N$ on $M\setminus\{v\}$, locally near $v$ in the splitting of variables 
$(r,x)\in\R_+\times X$ of the form
\begin{equation} \label{Fuchsint}
A:= r^{-\mu}\sum_{j=0}^\mu a_j(r)\left( -r\frac{\partial}{\partial r}\right)^j
\end{equation}
with coefficients $a_j\in C^\infty(\Rb_+,\textup{Diff}^{\mu-j}(X))$ (here $\textup{Diff}^\nu(\cdot)$ 
denotes the space of all differential operators of order $\nu$ on the manifold in parentheses, 
with smooth coefficients). Observe that when we consider a Riemannian metric on 
$\R_+\times X:=X^\wedge$ of the form $dr^2+r^2g_X$, where $g_X$ is a Riemannian metric on $X$, 
then the associated Laplace-Beltrami operator is just of the form \eqref{Fuchsint} for $\mu=2$. 
For such operators we have the homogeneous principal symbol $\sigma_\psi(A)\in C^\infty (T^*(M
\setminus\{v\})\setminus 0)$, and locally near $v$ in the variables $(r,x)$ 
with covariables $(\rho,\xi)$ the function
\[ \td\sigma_\psi (A)(r,x ,\rho ,\xi ):= r^\mu\sigma_\psi (A)(r,x,r^{-1}\rho ,\xi ) \]
which is smooth up to $r=0$. If a symbol (or an operator function) contains $r$ and $\rho$ in the combination $r\rho$ we speak of degeneracy of Fuchs type.\\
It is interesting to ask the nature of an operator algebra that contains Fuchs type differential operators of the from \eqref{Fuchsint} on $X^\Delta$, together with the parametrices of elliptic elements. An analogous problem is meaningful on $M$. Answers may be found in \cite{Schu2}, including the tools of the resulting so-called cone algebra. As noted above the ellipticity close to the tip $r=0$ is connected with a second symbolic structure, namely, the conormal symbol
\begin{equation} \label{co.}
\sigma_\textup{c}(A)(w):=\sum_{j=0}^\mu a_j(0)w^j:H^s(X)\to H^{s-\mu}(X) 
\end{equation}
which is a family of operators, depending on $w\in\Gamma_{\frac{n+1}{2}-\gamma}$, $\Gamma_\beta:= \{w\in\C:\reb w=\beta\}$, $n=\dim X$. Here $H^s(X)$ are the standard Sobolev spaces of smoothness $s\in\R$ on $X$. Ellipticity of $A$ with respect to a weight $\gamma\in\R$ means that \eqref{co.} is a family of isomorphisms for all $w\in\Gamma_{\frac{n+1}{2}-\gamma}$.\\
The ellipticity on the infinite cone $X^\Delta$ refers to a further principal symbolic structure, to be observed when $r\to\infty$. The behaviour in that respect is not symmetric under the substitution $r\to r^{-1}$. The present axiomatic approch will refer to ``abstract'' corners represented by $r\to 0$. The considerations are  based on specific insight on families of reductions of orders in given scales of spaces (in the simplest case $H^s(X), s\in\R$, when the corner is a conical sigularity). In order to motivate  our general constructions we briefly recall the form of corner operators of second generation. \\ 
\indent First, a differential operator on an open stretched wedge $\R_+\times X\times\Omega\ni(r,x,y)$, $\Omega\subseteq\R^q$ open, is called edge-degenerate, if it has the form
\begin{equation} \label{eeddeq}
A=r^{-\mu}\sum_{j+|\alpha|\leq\mu}a_{j\alpha}(r,y)\left(-r\frac{\partial}{\partial r}\right)^j (rD_y)^\alpha ,
\end{equation}
$a_{j\alpha}\in C^\infty(\Rb_+\times\Omega,\textup{Diff}^{\mu-(j+|\alpha|)}(X))$. Observe that \eqref{eeddeq} can be written in the form $A=r^{-\mu}\op_{r,y}(p)$ for an operator-valued symbol $p$ of the form $p(r,y,\rho,\eta)=\td p(r,y,r\rho,r\eta)$ and $\td p(r,y,\td\rho,\td\eta)\in C^\infty (\Rb_+\times\Omega,L_\clw^\mu(X;\R^{1+q}_{\td\rho,\td\eta}))$, 
\[\op_{r,y}(p)u(r,y)= \iint e^{i(r-r')\rho+i(y-y')\eta} p(r,y,\rho,\eta) u(r',y') dr'dy' \dbar\rho\dbar\eta . \]
Here $L_\clw^\mu(X;\R_\lambda^l)$ means the space of classical parameter-dependent pseudo-differen\-tial operators on $X$ of order $\mu$, with parameter $\lambda\in\R^l$, that is, locally on $X$ the operators are given in terms of amplitude functions $a(x,\xi,\lambda)$, where $(\xi,\lambda)$ is treated as an $(n+l)$-dimensional covariable, and we have $L^{-\infty}(X;\R^l):=\Sch(\R^l,L^{-\infty}(X))$ with $L^{-\infty}(X)$ being the (Fr\'echet) space of smoothing operators on $X$.\\ 
Let $\Diff_\textup{deg}^\mu(M)$ for a manifold $M$ with edge $Y$ denote the space of all differential operators on $M\setminus Y$ of order $\mu$ that are locally near $Y$ in the splitting of variables $(r,x,y)\in\R_+\times X\times\Omega$ of the form \eqref{eeddeq}. If we replace in the definition the edge covariable $\eta$ by $(\eta,\lambda)\in\R^{q+l}$ ($q=\dim Y$) we obtain parameter-dependent  families of operators in $\Diff_\textup{deg}^\mu(M)$. Similarly as \eqref{Fuchsint} an operator of the form
\[ A:=t^{-\mu}\sum_{j=0}^\mu a_j(t)\left(-t\frac{\partial}{\partial t}\right)^j \]
is called corner degenerate  if $a_j\in C^\infty(\Rb_+,\Diff_\textup{deg}^{\mu-j}(M))$, $j=0,1,\ldots, \mu$. The corner conormal symbol $\sigma_\textup{c}(A)(z)=\sum_{j=0}^\mu a_j(0)z^j$, $z\in \Gamma_{\frac{\dim M +1}{2}-\delta}$ for a corner weight $\delta\in\R$, is just a parameter-dependent family in $\Diff_\textup{deg}^\mu(M)$ with parameter $\imb z$ on the indicated weight line. The program to study such operators close to the tip $t\to 0$ (see \cite{Calv2}, \cite{Haru11}) is just a concrete realisation of the present theory.\\
\indent This paper is organised as follows. In Chapter \ref{ElofCo} we introduce spaces of symbols based on families of reductions of orders in given scales of (analogues of Sobolev) spaces.\\
Chapter \ref{Oprefex} is devoted to the specific effects of an axiomatic calculus near the tip of the corner. The corner axis is represented by a real axis $\R\ni r$, and the operators take values in vector-valued analogues of Sobolev spaces in $r$.\\
\indent As indicated above, our results are designed  as a step of a larger concept of abstract edge and corner theories, organised in an iterative manner. The full calculus employs the one for $r\to\infty$ as a counterpart of our Mellin operators on $\R_+$ near $r=0$. However, the continuation of the calculus in that sense needs more space than available in the present note.

\section{Symbols associated with order reductions}
\label{ElofCo}
\subsection{Scales and order reducing families}
\label{OrdRed}
Let $\mathfrak{E}$ denote the set of all families $\E=(E^s)_{s\in\R}$ of Hilbert spaces with continuous embeddings $E^{s'}\hookrightarrow E^s$, $s'\geq s$, so that $E^\infty:=\bigcap _{s\in\R}E^s$ is dense in every $E^s,s\in\R$ and that there is a dual scale $\E^*=(E^{*s})_{s\in\R}$ with a non-degenerate sesquilinear pairing $(.,.)_0:E^0\times E^{*0}\to\C$, such that $(.,.)_0:E^\infty\times E^{*\infty}\to\C$, extends to a non-degenerate sesquilinear pairing
\[ E^s\times E^{*-s}\to\C \]
for every $s\in\R$, where $\sup_ {f\in E^{*-s}\setminus \{ 0 \} }\frac{|(u,f)_0|}{\|f\|_{E^{*-s}}}$ and $\sup_{g\in E^s\setminus \{ 0 \} }\frac{|(g,v)_0|}{\|g\|_{E^s}}$ are equivalent norms in the spaces $E^s$ and $E^{*-s}$, respectively; moreover, if $\E=(E^s)_{s\in\R}$, $\widetilde\E=(\widetilde E^s)_{s\in\R}$ are two scales in consideration and 
$a\in\Li^\mu(\E,\widetilde\E):=\bigcap_{s\in\R}\Li(E^s,\widetilde E^{s-\mu})$, for some $\mu\in\R$, then
\[ \sup_{s\in [s',s'']} \|a\|_{s,s-\mu} < \infty \]
for every $s'\leq s''$; here $\|.\|_{s,\td s}:=\|.\|_{\Li(E^s,\widetilde E^{\td s})}$. Later on, in the case $s=\td s= 0$ we often write $\|.\|:=\|.\|_{0,0}$.\\
Let us say that a scale $\E\in\mathfrak{E}$ is said to have the compact embedding property, if the embeddings 
$E^{s'}\hookrightarrow E^s$ are compact when $s'>s$.
\begin{Rem}
Every $a\in\Li^\mu(\E,\widetilde\E)$ has a formal adjoint $a^*\in\Li^\mu(\widetilde\E^*,\E^*)$, obtained by $(au,v)_0=(u,a^*v)_0$ for all $u\in E^\infty,v\in \widetilde E^{*\infty}$.
\end{Rem}
\begin{Rem}
The space $\Li^\mu(\E,\widetilde\E)$ is Fr\'{e}chet in a natural way for every $\mu\in\R$. 
\end{Rem}
\begin{Def} \label{ordred}
A system $(b^\mu(\eta))_{\mu\in\R}$ of operator functions $b^\mu (\eta)\in C^\infty (\R^{q},\Li^\mu 
(\E,\E ))$ is called an order reducing family of the scale $\E$, if $b^\mu(\eta):E^s\to E^{s-\mu}$ is an isomorphism for every $s,\mu\in\R$, $\eta\in\R^q$, $b^0(\eta)=\textup{id}$ for every $\eta\in\R^q$, and
\begin{enumerate}
\item[\textup{(i)}]$D_{\eta}^\beta b^\mu(\eta)\in C^\infty (\R^{q},\Li^{\mu
-|\beta|}(\E,\E))$ for every $\beta\in\N^{q}$; 
\item[\textup{(ii)}] for every $s\in\R,\beta\in\N^q$ we have 
\[ \max_{|\beta|\leq k}\sup_{\eta\in\R^q \atop s\in [s',s''] } \|b^{s-\mu+|\beta|}(\eta)\{D_\eta^\beta b^\mu(\eta )\}b^{-s}(\eta )\|_{0,0}< \infty \]
for all $k\in\N$, and for all real $s'\leq s''$.
\item[\textup{(iii)}] for every $\mu,\nu\in\R$, $\nu\geq\mu$, we have
\[ \sup_{s\in [s',s'']}\|b^\mu(\eta)\|_{s,s-\nu}\leq c\la\eta\ra^B \]
for all $\eta\in\R^q$ and $s'\leq s''$ with constants $c(\mu,\nu,s), B(\mu,\nu,s)>0$, uniformly bounded in compact $s$-intervals and compact $\mu,\nu$-intervals for $\nu\geq\mu$; moreover, for every $\mu\leq 0$ we have 
\[ \|b^\mu(\eta)\|_{0,0}\leq c\la\eta\ra^\mu \]
for all $\eta\in\R^q$ with constants $c>0$, uniformly bounded in compact $\mu$-intervals, $\mu\leq 0$. 
\end{enumerate}
\end{Def}
Clearly the operators $b^\mu$ in (iii) for $\nu\geq\mu$ or $\mu\leq 0$, are composed with a corresponding 
embedding operator. \\
In addition we require that the operator families $\left(b^\mu(\eta)\right)^{-1}$ are equivalent to $b^{-\mu}(\eta)$, according to the following notation. Another order reducing family $(b_1^\mu(\eta))_{\mu\in\R},\eta\in\R^q$, 
in the scale $\E$ is said to be equivalent to $(b^\mu(\eta))_{\mu\in\R}$, if for every $s\in\R,
\beta\in\N^q$, there are constants $c=c(\beta,s)$ such that 
\[ \|b_1^{s-\mu+|\beta|}(\eta)\{D_\eta^\beta b^\mu(\eta)\}b_1^{-s}(\eta)\|_{0,0}\leq c , \]
\[ \|b^{s-\mu+|\beta|}(\eta)\{D_\eta^\beta b_1^\mu(\eta)\}b^{-s}(\eta)\|_{0,0}\leq c , \]
for all $\eta\in\R^q$, uniformly in $s\in [s',s'']$ for every $s'\leq s''$.
\begin{Rem}
\label{parrem}
Parameter-dependent theories of operators are common in many concrete contexts. For instance, if $\Omega$ is an \lr open\rr\ $C^\infty$ manifold, there is the space $L_\clw^\mu(\Omega,\R^q)$ of parameter-dependent pseudo-differential operators on $\Omega$ of order $\mu\in\R$, with parameter $\eta\in\R^q$, where the local amplitude functions $a(x,\xi,\eta)$ are classical symbols in $(\xi,\eta)\in\R^{n+q}$, treated as covariables, $n=\dim\Omega$, while $L^{-\infty}(\Omega,\R^q)$ is the space of Schwartz functions in $\eta\in\R^q$ with values in $L^{-\infty}(\Omega)$, the space of smoothing operators on $\Omega$. Later on we will also consider specific examples with more control on the dependence on $\eta$, namely, when $\Omega=M\setminus\{v\}$ for a manifold $M$ with conical singularity $v$.
\end{Rem}
\begin{Exam} \label{ExamX}
Let $X$ be a closed compact $C^\infty$ manifold, $E^s:=H^s(X),s\in\R$, the scale of classical 
Sobolev spaces on $X$ and $b^\mu(\eta)\in L_\clw^\mu (X;\R^q_\eta)$ a parameter-dependent elliptic family that induces isomorphisms $b^\mu(\eta):H^s(X)\to H^{s-\mu}(X)$ for all $s\in\R$. Then for 
$\nu\geq\mu$ we have
\[ \|b^\mu(\eta)\|_{\Li(H^s(X),H^{s-\nu}(X))}\leq c\la\eta\ra^{\pi(\mu,\nu)} \]
for all $\eta\in\R^q$, uniformly in $s\in [s',s'']$ for arbitrary $s',s''$, as well as in compact $\mu$- and $\nu$-intervals for $\nu\geq\mu$, where
\begin{equation} 
\pi(\mu,\nu):=\max(\mu,\mu-\nu)
\end{equation}
with a constant $c=c(\mu,\nu,s',s'')>0$. Observe that $\sup_{\xi\in\R^p} \frac{\la\xi,\eta\ra^\mu} {\la\xi\ra^\nu}\leq\la\eta\ra^{\pi(\mu,\nu)}$
for all $\eta\in\R^q$.
\end{Exam}
\begin{Rem}
Let $b^s(\td\tau,\td\eta)\in L^\mu_\clw(X;\R_{\td\tau,\td\eta}^{1+q})$ be an order reducing family as in the above example, now with the parameter $(\td\tau,\td\eta)\in\R^{1+q}$ rather than $\eta$, and of order $s\in\R$. Then, setting $b^s(t,\tau,\eta):=b^s(t\tau,t\eta)$ the expression
\[ \Big\{ \int \| [t]^{-s}\op_t(b^s)(\eta^1)u\|^2_{L^2(X)}dt\Big\}^\oneh \]
for $\eta^1\in\R^q\setminus\{0\}$, $|\eta^1|$ sufficiently large, is a norm on the space $\Sch(\R,C^\infty(X))$. Let $H^s_\cone(\R\times X)$ denote the completion of $\Sch(\R,C^\infty(X))$ in this norm. Observe that this space is independent of the choice of $\eta^1$, $|\eta^1|$ sufficiently large. For reference below we also form weighted variants $H^{s;g}_\cone(\R\times X):=\la t\ra^{-g}H^s_\cone (\R\times X)$, $g\in\R$, and set 
\begin{equation} \label{equtcon}
H^{s;g}_\cone(\R_+\times X):=H^{s;g}_\cone(\R\times X)|_{\R_+\times X}.
\end{equation}
As is known, cf. \textup{\cite{Kapa10}}, the spaces $H^{s;g}_\cone(\R\times X)$ are weighted Sobolev spaces in the calculus of pseudo-differential operators on $\R_+\times X$ with $|t|\to\infty$ being interpreted as a conical exit to infinity.
\end{Rem}
Another feature of order reducing families, known, for instance, in the case of the above example, is that when $U\subseteq\R^p$ is an open set and $m(y)\in C^\infty(U)$ a strictly positive function, $m(y)\geq c$ for $c>0$ and for all $y\in U$, the family $b_1^s(y,\eta):=b^s(m(y)\eta)$, $s\in\R$, is order reducing in the sense of Definition \ref{ordred} and equivalent to $b(\eta)$ for every $y\in U$, uniformly in $y\in K$ for any compact subset $K\subset U$. A natural requirement is that when $m>0$ is a parameter, there is a constant $M=M(s',s'')>0$ such that
\begin{equation} \label{steeq} 
\| b^s(\eta)b^{-s}(m\eta)\|_{0,0}\leq c\max(m,m^{-1})^M
\end{equation}
 for every $s\in[s',s'']$, $m\in\R_+$, and $\eta\in\R^q$.\\
We now turn to another example of an order reducing family, motivated by the calculus of pseudo-differential operators on a manifold with edge (here in ``abstract'' form), where all the above requirements are satisfied, including the latter one.
\begin{Def} \label{DefHkap}
\begin{enumerate}
\item[\textup{(i)}] If $H$ is a Hilbert space and $\kappa:=\{\kappa_\lambda\}_{\lambda\in\R_+}$ a group of isomorphisms $\kappa_\lambda:H\to H$, such that $\lambda\to\kappa_\lambda h$ defines a continuous function $\R_+\to H$ for every $h\in H$, and $\kappa_\lambda\kappa_\rho=\kappa_{\lambda\rho}$ for $\lambda,\rho\in\R$, we call $\kappa$ a group action on $H$.
\item[\textup{(ii)}] Let $\Hi=(H^s)_{s\in\R}\in\mathfrak{E}$ and assume that $H^0$ is endowed with a group action 
$\kappa=\{\kappa_\lambda\}_{\lambda\in\R_+}$ that restricts \lr for $s>0$\rr\ or extends \lr 
for $s<0$\rr\ to a group action on $H^s$ for every $s\in\R$. In addition, we assume that $\kappa$ 
is a unitary group action on $H^0$. We then say that $\Hi$ is endowed with a group action. 
\end{enumerate}
\end{Def}
If $H$ and $\kappa$ are as in Definition \ref{DefHkap} (i), it is known that there are constants $c,M>0$, such that 
\begin{equation} \label{kappMc}
\| \kappa_\lambda\|_{\Li(H)}\leq c\max (\lambda,\lambda^{-1})^M
\end{equation}
for all $\lambda\in\R_+$.\\
Let $\We^s(\R^q,H)$ denote the completion of $\Sch(\R^q,H)$ with respect to the norm
\[ \|u\|_{\We^s(\R^q,H)}:=\Big\{\int \la\eta\ra^{2s}\|\kappa^{-1}_{\la\eta\ra}\hat{u}(\eta)\|_H^2
d\eta\Big\}^\oneh ; \]
$\hat{u}(\eta)=\Fu_{y\to\eta}u(\eta)$ is the Fourier transform in $\R^q$. The space $\We^s(\R^q,H)$ 
will be referred to as edge space on $\R^q$ of smoothness $s\in\R$ (modelled on $H$). 
Given a scale $\Hi=(H^s)_{s\in\R}\in\mathfrak{E}$ with group action we have the edge spaces
\[ W^s:=\We^s(\R^q,H^s),\ s\in\R. \]
If necessary we also write $\We^s(\R^q,H^s)_\kappa$.   
The spaces form again a scale $\We:=(W^s)_{s\in\R}\in\mathfrak{E}$.\\
For purposes below we now formulate a class of operator-valued symbols 
\begin{equation} \label{neu5}
S^\mu (U\times\R^q ;H,\widetilde{H})_{\kappa,\td\kappa }
\end{equation}
for open $U\subseteq\R^p$ and Hilbert spaces $H$ and $\widetilde{H}$, endowed with group actions $\kappa=\{
\kappa_\lambda\}_{\lambda\in\R_+}, \td\kappa=\{\td\kappa_\lambda\}_{\lambda\in\R_+}$, respectively, as 
follows. The space \eqref{neu5} is defined to be the set of all $a(y,\eta)\in C^\infty(U\times\R^q,
\Li(H,\widetilde{H}))$ such that
\begin{equation} \label{neu2,5}
\sup_{(y,\eta )\in K\times\R^q }\la\eta\ra^{-\mu+|\beta|}\|\td\kappa^{-1}_{\la\eta\ra}\{D_y^\alpha
D_\eta^\beta a(y,\eta )\} \kappa_{\la\eta\ra}\|_{\Li(H,\td H)} < \infty
\end{equation}
for every $K\Subset U,\alpha\in\N^p,\beta\in\N^q$.
\begin{Rem} \label{rem7}
Analogous symbols can also be defined in the case when $\widetilde H$ is a Fr\'{e}chet space with group action, i.e., $\widetilde H$ is written as a projective limit of Hilbert spaces $\widetilde H_j, j\in\N$, with continuous embeddings $\widetilde H_j\hookrightarrow\widetilde H_0$, where the group action on $\widetilde H_0$ restricts to group actions on $\widetilde H_j$ for every $j$. Then $S^\mu(U\times\R^q;H,\widetilde H):= \varinjlim_{j\in\N} S^\mu(U\times\R^q;H,\widetilde H_j)$.
\end{Rem}
Consider an operator function $p(\xi,\eta)\in C^\infty (\R^{p+q}_{\xi,\eta},\Li^\mu(\Hi,\Hi))$
that represents a symbol
\[ p(\xi,\eta)\in S^\mu(\R^{p+q}_{\xi,\eta};H^s,H^{s-\mu})_{\kappa,\kappa} \]
for every $s\in\R$, such that $p(\xi,\eta):H^s\to H^{s-\mu}$ is a family of isomorphisms for all $s\in\R$, and the inverses $p^{-1}(\xi,\eta)$ represent a symbol 
\[ p^{-1}(\xi,\eta)\in S^{-\mu}(\R^{p+q}_{\xi,\eta};H^s,H^{s+\mu})_{\kappa,\kappa} \]
for every $s\in\R$. Then $b^\mu(\eta):=\op_x(p)(\eta)$ is a family of isomorphisms
\[ b^\mu(\eta):W^s\to W^{s-\mu},\ \eta\in\R^q , \]
with the inverses $b^{-\mu}(\eta):=\op_x(p^{-1})(\eta)$.
\begin{Prop} \label{nn9}
\begin{enumerate}
\item[\textup{(i)}] We have
\begin{equation} \label{731.a0802.eq}
\|b^\mu(\eta)\|_{\Li(W^0,W^0)}\leq c\la\eta\ra^\mu 
\end{equation}
for every $\mu\leq 0$, with a constant $c(\mu)>0$.
\item[\textup{(ii)}] For every $s,\mu,\nu\in\R$, $\nu\geq\mu$, we have
\begin{equation}
\|b^\mu(\eta)\|_{\Li(W^s,W^{s-\nu})}\leq c\la\eta\ra^{\pi(\mu,\nu)+M(s)+M(s-\mu)}
\end{equation}
for all $\eta\in\R^q$, with a constant $c(\mu,s)>0$, and $M(s)\geq 0$ defined by
\[ \|\kappa_\lambda\|_{\Li(H^s,H^s)}\leq c\lambda^{M(s)}\ \textup{for all}\ \lambda\geq 1. \]
\end{enumerate}
\end{Prop}  
\proof 
(i) Let us check the estimate \eqref{731.a0802.eq}. 
For the computations we denote by $j : H^{- \mu}  \hookrightarrow
 H^0$ the embedding
operator. We have for $u \in W^0$
\begin{align*}
\| b^\mu(\eta) u \|_{W^0}^2 = & \int \| j p(\xi, \eta) (\Fu u)(\xi) \|_{H^0}^2 d \xi  \\
= & \int \| \kappa_{ \langle \xi, \eta \rangle }^{-1} j 
       \kappa_{
       \langle \xi, \eta \rangle } \kappa_{ \langle 
       \xi, \eta \rangle }^{-1}    
       p(\xi, \eta) \kappa_{ \langle \xi, \eta
       \rangle }  
       \kappa_{ \langle \xi, \eta \rangle }^{-1} (\Fu u)(\xi) \|_{H^0}^2 d \xi \\
\leq & \int \| \kappa_{ \langle \xi, \eta \rangle }^{-1} j \kappa_{
        \langle \xi, \eta \rangle } \|^2_{ {\Li}(H^{- \mu},H^0)}
     \| \kappa_{ \langle \xi, \eta \rangle }^{-1} 
        p(\xi, \eta) \kappa_{ \langle \xi, \eta 
	\rangle } \kappa_{ \langle \xi,
     \eta \rangle }^{-1} (\Fu u)(\xi) \|_{H^{- \mu}}^2 d \xi \\
\leq & c \int \| \kappa_{ \langle \xi, \eta \rangle }^{-1} p(\xi,
       \eta) \kappa_{ \langle \xi, \eta \rangle } \|^2_{ 
       {\Li}(H^0, H^{- \mu})} \| \kappa_{ \langle \xi, \eta
       \rangle }^{-1} (\Fu u)(\xi) \|_{H^0}^2 d \xi   \\
\leq & c \sup_{\xi \in \R^p} \langle \xi, \eta \rangle^{2\mu} \|
     u \|_{W^0}^2. \\
\end{align*}
Thus $\| b^\mu(\eta) \|_{ {\mathcal{L}}(W^0, W^0)} \leq c \sup_{\xi \in
      \R^p} \langle \xi, \eta \rangle^\mu  \leq c \langle \eta
      \rangle^\mu$, since $\mu \leq 0$.

(ii) Let $j : H^{s-\mu} \hookrightarrow H^{s-\nu}$ denote the
canonical embedding.  
For every fixed $s \in \R$ we have
\begin{multline*}
\| b^\mu(\eta) u \|_{W^{s-\nu}}^2 \ =  \int \langle \xi
      \rangle^{2(s-\nu)} \| \kappa_{ \langle \xi
      \rangle}^{-1}
      j p(\xi, \eta)(\Fu_{x \to \xi} u)(\xi)
      \|_{H^{s- \nu}}^2 d \xi  \\
   =    \int \langle \xi \rangle^{2(s-\nu)} \| \kappa_{
         \langle \xi \rangle }^{-1} j p(\xi, \eta)
	 \kappa_{ \langle \xi \rangle }
         \langle \xi \rangle^{-s} \langle \xi \rangle^s 
	 \kappa_{ \langle \xi \rangle }^{-1}
	 (\Fu_{x \to \xi} u)(\xi) \|_{H^{s-\nu}}^2 d 
	     \xi   \\
    =   \sup_{\xi \in \R^p} \langle \xi \rangle^{-2 \nu} \|
         \kappa_{ \langle \xi \rangle }^{-1} j 
	 p(\xi, \eta) \kappa_{ \langle \xi \rangle } 
	 \|_ { {\Li}(H^s, H^{s-\nu})}^2    
   \int \langle \xi \rangle^{2s} \| \kappa_{ \langle \xi
          \rangle }^{-1} \Fu_{x \to \xi} u(\xi) 
	  \|_{H^s}^2 d \xi
\end{multline*}
We have
\begin{multline*}
\| \kappa_{ \langle \xi \rangle }^{-1} 
     \bigl(j p(\xi, \eta)\bigr) \kappa_{ \langle \xi
     \rangle } \|_{ {\mathcal{L}}(H^s, H^{s-\nu})}   \\
  \leq   \| \kappa_{ \langle \xi \rangle }^{-1} j
            \kappa_{ \langle \xi \rangle } \|_{ 
	  {\mathcal{L}} (H^{s-\mu}, H^{s- \nu})} 
   \| \kappa_{ \langle \xi \rangle }^{-1} p(\xi,
        \eta) \kappa_{ \langle \xi \rangle }     %^{-1} 
	 \|_{ {\mathcal{L}}(H^s, H^{s-\mu})} \\
\leq  c \| \kappa_{ \langle \xi \rangle}^{-1} p 
         (\xi, \eta)
         \kappa_{ \langle \xi \rangle }
	 \|_{ {\mathcal{L}}(H^s, H^{s- \mu})}	  
\end{multline*}
with a constant $c> 0$.	 

We employed here that $\| \kappa_{ \langle \xi 
\rangle }^{-1} j
\kappa_{ \langle \xi \rangle } \|_{ {\mathcal{L}}(H^{s-\mu}, 
H^{s-\nu})} \leq c$
for all $\xi \in \R^p$. Moreover, we have
\begin{multline*}
\| \kappa_{ \langle \xi \rangle }^{-1}  p(\xi, \eta)
     \kappa_{ \langle \xi \rangle } \|_{ {\mathcal{L}}(H^s, H^{s-
     \mu})}   \\
\leq  \| \kappa_{ \langle \xi \rangle }^{-1} \kappa_{ \langle
     \xi,\eta \rangle } \|_{ {\mathcal{L}}(H^{s-\mu}, H^{s-\mu})}
     \| \kappa_{ \langle \xi, \eta \rangle }^{-1} p 
     (\xi,\eta)  \kappa_{ \langle \xi, \eta \rangle} 
     \|_{ {\mathcal{L}}(H^s, H^{s-\mu})}
     \| \kappa_{ \langle \xi, \eta \rangle}^{-1} \kappa_{
     \langle \xi \rangle} \|_{ {\mathcal{L}}(H^s, H^s)}    \\
\leq  c \langle \xi, \eta \rangle^\mu \| \kappa_{\langle \xi,
     \eta \rangle  \langle \xi \rangle^{-1}} \|_{ {\mathcal{L}}(H^{s-
     \mu}, H^{s- \mu})} \| \kappa_{ \langle \xi, \eta
     \rangle^{-1}\langle \xi \rangle } \|_{ 
        {\mathcal{L}}(H^s, H^s)}   \\
\leq  c \langle \xi, \eta\rangle^\mu 
        \bigl( \frac{\langle \xi,
       \eta \rangle}{\langle \xi \rangle} \bigr)^{M(s-\mu) +
       M(s)}.
\end{multline*}           

As usual, $c > 0$ denotes different constants (they may also
depend on $s$); the numbers $M(s)$, $s \in \R$, are determined
by the estimates
\[  \| \kappa_\lambda \|_{ {\mathcal{L}}(H^s, H^s)} \leq c
       \lambda^{M(s)} \ \textup{for all $\lambda \geq 1$}.    \]
We obtain altogether that       
\[  \| b^\mu(\eta) \|_{ {\mathcal{L}}(W^s,W^{s-\nu})} \leq
       c \sup_{\xi \in \R^n} \frac{\langle \xi, \eta
       \rangle^\mu}{\langle \xi \rangle^\nu} \bigl(
       \frac{\langle \xi, \eta\rangle}{\langle \xi \rangle}
       \bigr)^{M(s-\mu)+M(s)}
       \leq c \langle \eta \rangle^{\pi(\mu, \nu) +
       M(s-\mu)+M(s)}.      \]
\qed \\
It can be proved that the operators in Proposition \ref{nn9} also have the uniformity properties with respect to $s,\mu,\nu$ in compact sets, imposed in Definition \ref{ordred}.

\subsection{Symbols based on order reductions}
\label{opvalsymb}
We now turn to operator valued symbols, referring to scales
\[ \E=(E^s)_{s\in\R}, \widetilde{\E}=(\widetilde{E}^s)_{s\in\R}\in\mathfrak{E}. \]
For purposes below we slightly generalise the concept of order reducing families by replacing the parameter space $\R^q\ni\eta$ by $\He\ni\eta$, where
\begin{equation} \label{hpar}
\He:=\{\eta=(\eta',\eta'')\in\R^{q'+q''}:q=q'+q'', \eta'' \not=0\}.
\end{equation}
In other words for every $\mu\in\R$ we fix order-reducing families $b^\mu(\eta)$ and $\td b^\mu(\eta)$ in the scales $\E$ and $\widetilde{\E}$, respectively, where $\eta$ varies over $\He$, and the properties of Definition \ref{ordred} are required for all $\eta\in \He$. In many cases we may admit the case $\He=\R^q$ as well.
\begin{Def} \label{732.symb1901.de}
By $S^\mu(U\times\He;\E,\widetilde\E)$ for open $U\subseteq \R^p,\mu\in\R$, we denote the set of all $a(y,\eta )\in C^\infty ( U\times\He ,\Li^\mu (\E,\widetilde\E ))$ such that
\begin{equation} \label{732.est1901.eq} 
D^\alpha_y D^\beta_\eta a(y,\eta )\in C^\infty (U\times\He ,\Li^{\mu-|\beta|}(\E,\widetilde{\E} )), 
\end{equation}
and for every $s\in\R$ we have 
\begin{equation} \label{opvasy} 
\max_{|\alpha|+|\beta|\leq k} \sup_{y\in K , \eta\in\He, \eta\geq h \atop s\in [s',s''] } \| \td b^{s-\mu +|\beta |}(\eta )
\{ D^\alpha_y D^\beta_\eta a(y,\eta )\} b^{-s}(\eta )\|_{0,0}
\end{equation}
is finite for all $K\Subset U$, $k\in\N, h>0$.
\end{Def}
Let $S^\mu(\He;\E,\widetilde{\E})$ denote the subspace of all elements of $S^\mu(U\times\He;\E,\widetilde{\E})$ 
that are independent of $y$.\\
Observe that when $(b^\mu(\eta))_{\mu\in\R}$ is an order reducing family parametrised by $\eta\in\He$ then we have
\begin{equation} \label{ordsmb} 
b^\mu (\eta )\in S^\mu (\He ; \E ,\E )
\end{equation}
for every $\mu\in\R$.
\begin{Rem} \label{Frechnat}
The space $S^\mu(U\times\He;\E,\widetilde\E)$ is Fr\'{e}chet with the 
semi-norms 
\begin{equation} \label{semism}
a\to\max_{|\alpha|+|\beta|\leq k}\sup_{(y,\eta)\in K\times\He, |\eta|\geq h \atop s\in [s',s''] }\|\td b^{s-\mu+|\beta|}(\eta)\{D_y^\alpha D_\eta^\beta a(y,\eta)\} b^{-s}(\eta)\|_{0,0}
\end{equation}
parametrised by $K\Subset U$, $s\in\Z$, $\alpha\in\N^p$, $\beta\in\N^{q}$, $h>0$, which are the best constants in the estimates \eqref{opvasy}. We then have
\[ S^\mu (U\times\He ;\E ,\widetilde\E )= C^\infty (U,S^\mu (\He ;\E,\widetilde\E )) = C^\infty 
(U)\hat\otimes_\pi S^\mu (\He ;\E ,\widetilde\E) . \]
\end{Rem} 
We will also employ other variants of such symbols, for instance, when $\Omega\subseteq\R^m$ is an open set, 
\[ S^\mu (\Rb_+\times\Omega\times\He;\E,\widetilde\E):= C^\infty(\Rb_+\times\Omega, S^\mu (\He;\E,\widetilde\E)). \]\\
\indent In order to emphasise the similarity of our considerations for $\He$ with the case $\He=\R^q$ we often write again $\R^q$ and later on tacitly use the corresponding results for $\He$ in general.
\begin{Rem} \label{Rem7315}
Let $a(y,\eta)\in S^\mu(U\times\R^q)$ be a polynomial in $\eta$ of order $\mu$ and $\E=(E^s)_ {s\in\R}$ a scale and identify $D_y^\alpha D_\eta^\beta a(y,\eta)$ with $\big(D_y^\alpha D_\eta^\beta a(y,\eta)\big)\iota$ with the embedding $\iota: E^s\to E^{s-\mu+|\beta|}$. Then we have
\begin{align*}
\| b^{s-\mu+|\beta|} & (\eta)\big(D_y^\alpha D_\eta^\beta a(y,\eta)\big)b^{-s}(\eta)\|_ {0,0} \\ & \leq |D_y^\alpha D_\eta^\beta a(y,\eta)| \|b^{-\mu+|\beta|}(\eta)\|_ {0,0} \leq c \la\eta\ra^{\mu-|\beta|}\la\eta\ra^{-\mu+|\beta|}=c
\end{align*}
for all $\beta\in\N^q$, $|\beta|\leq\mu$, $y\in K\Subset U$ \lr see Definition \textup{\ref{ordred}} \textup{(iii)}\rr. Thus $a(y,\eta)$ is canonically identified with an element of $S^\mu(U\times\R^q ;\E,\E)$.
\end{Rem}
\begin{Prop} \label{fs-}
We have 
\[ S^{-\infty}(U\times\R^q;\E,\widetilde\E):=\bigcap_{\mu\in\R}S^\mu(U\times\R^q;\E,\widetilde\E) = C^\infty(U,\Sch(\R^q,\Li^{-\infty}(\E,\widetilde\E))). \]
\end{Prop}
\proof
Let us show the assertion for $y$-independent symbols; the $y$-dependent case is then straightforward. For notational convenience we set $\widetilde\E=\E$; the general case is analogous.  First let $a(\eta)\in S^{-\infty}(\R^q;\E,\E)$, which means that $a(\eta)\in C^\infty(\R^q,\Li^{-\infty}(\E,\E))$ and 
\begin{equation} \label{equat1}
\|b^{s+N}(\eta)\{D_\eta^\beta a(\eta)\}b^{-s}(\eta)\|_{0,0}< c
\end{equation}
for all $s\in\R$, $N\in\N$, $\beta\in\N^q$ and show that
\begin{equation} \label{equat2}
\sup_{\eta\in\R^q}\|\la\eta\ra^M D_\eta^\beta a(\eta)\|_{s,t}<\infty
\end{equation}
for every $s,t\in\R$, $M\in\N$, $\beta\in\N^q$. To estimate \eqref{equat2} it is enough to assume $t>0$. We have
\begin{equation} \label{estimate}
\|\la\eta\ra^M D_\eta^\beta a(\eta)\|_{s,t}=\|b^{-kt}(\eta)b^{kt}(\eta)\la\eta\ra^M D_\eta^\beta a(\eta)b^{-s}(\eta)b^s(\eta)\|_{s,t}
\end{equation}
for every $k\in\N,k\geq 1$, it is sufficient to show that the right hand side is uniformly bounded in $\eta\in\R^q$ for sufficiently large choice of $k$. The right hand side of \eqref{estimate} can be estimated by
\[\|b^{-t}(\eta)\|_{0,t}\|b^{(1-k)t}(\eta)\|_{0,0}\|b^{kt}(\eta)D_\eta^\beta a(\eta) b^{-s}(\eta)\|_{0,0}\|b^s(\eta)\|_{s,0}. \]
Using $\|b^{kt}(\eta) D_\eta^\beta a(\eta)b^{-s}(\eta)\|_{0,0}\leq c$, which is true by assumption and the estimates
\[ \|b^s(\eta)\|_{s,0}\leq c\la\eta\ra^B, \ \ \|b^{-t}(\eta)\|_{0,t}\leq c\la\eta\ra^{B'}, \]
with different $B,B'\in\R$ and $\|b^{(1-k)t}(\eta)\|_{0,0} \leq c\la\eta\ra^{(1-k)t}$
(see Definition \ref{ordred} (iii)) we obtain altogether
\[ \|\la\eta\ra^M D_\eta^\beta a(\eta)\|_{s,t}\leq c\la\eta\ra^{M+B+B'+(1-k)t} \]
for some $c>0$. Choosing $k$ large enough it follows that the exponent on the right hand side is $<0$, i.e., we obtain uniform boundedness in $\eta\in\R^q$. \\
To show the reverse direction suppose that $a(\eta)$ satisfies \eqref{equat2}, and let $\beta\in
\N^q$, $M,s,t\in\R$ be arbitrary. We have
\begin{multline} \label{convdir}
\| b^t(\eta)D_\eta^\beta a(\eta)b^{-s}(\eta)\|_{0,0}\leq \\ \| b^t(\eta)\la\eta\ra^{-M}\|_{t,0}\|\la\eta\ra^{2M}D_\eta^\beta a(\eta)\|_{s,t}\|\la\eta\ra^{-M}b^{-s}(\eta)\| 
_{0,s}.
\end{multline}
Now using \eqref{equat2} and the estimates
\[ \|b^t(\eta)\la\eta\ra^{-M}\|_{t,0}\leq c \la\eta\ra^{A-M}, \ \ 
   \|\la\eta\ra^{-M}b^{-s}(\eta)\|_{0,s}\leq c \la\eta\ra^{A'-M}, \]
with constants $A,A'\in\R$, we obtain
\[ \| b^t(\eta)D_\eta^\beta a(\eta)b^{-s}(\eta)\|_{0,0}\leq c\la\eta\ra^{A+A'-2M}. \]
Choosing $M$ large enough we get uniform boundedness of \eqref{convdir} in $\eta\in\R^q$ which completes the proof.
\qed
\begin{Prop} \label{732.a2201.pr}
Let $a(y,\eta)\in S^\mu(U\times\R^q;\E,\widetilde\E)$ and $\mu\leq 0$. Then we have
\[ \| a(y,\eta )\|_{0,0}\leq c\la\eta\ra^\mu \]
for all $y\in K\Subset U,\eta\in\R^q$, with a constant $c=c(s,K)>0$.
\end{Prop}
\proof
For simplicity we consider the $y$-independent case. It is enough to show that $\| a(\eta ) u \|_{\widetilde{E}^0} \leq c \la\eta\ra^\mu \| u \|_{E^0}$ for all $u \in E^\infty$. Let $j : E^{-\mu} \to E^0$ denote the embedding operator. We then have  
\begin{align*}
\| a(\eta ) u \|_{ \widetilde{E}^0} = & 
\| a(\eta ) b^{-\mu }(\eta ) j b^\mu(\eta) u\|_{ \widetilde{E}^0} \\
\leq & \| a(\eta ) b^{-\mu }(\eta ) \|_{\Li (E^0,\widetilde{E}^0)} \| j b^\mu (\eta ) u \|_{E^0}
\leq c \la\eta\ra^\mu \| u\|_{E^0}.
\end{align*}
\qed
\begin{Prop} \label{s,s-nu}
A symbol $a(y,\eta)\in S^\mu(U\times\R^q;\E,\widetilde\E),\ \mu\in\R$, satisfies the estimates
\begin{equation} \label{symbest}
\|a (y,\eta )\|_{s,s-\nu}\leq c\la\eta\ra^A
\end{equation}
for every $\nu\geq\mu$, for every $y\in K\Subset U,\eta\in\R^q,s\in\R$, with constants $c=c
(s,\mu,\nu)>0,A=A(s,\mu,\nu,K)>0$ that are uniformly bounded when $s,\mu,\nu$ vary over compact sets, $\nu\geq\mu$.
\end{Prop}
\proof
For simplicity we consider again the $y$-independent
case. Let $j: \widetilde{E}^{s-\mu} \hookrightarrow
\widetilde{E}^{s-\nu}$ be the embedding operator. Then we
have
\begin{multline*}
\| a(\eta) \|_{s,s-\nu} =   \| j \tilde{b}^{-s+\mu}(\eta )
    \tilde{b}^{s-\mu}(\eta )a(\eta ) b^{-s}(\eta ) b^s(\eta )
     \|_{s,s-\nu}   \\
\leq  \| j \tilde{b}^{-s+\mu } (\eta ) \|_{0,s-\nu} \|
       \tilde{b}^{s- \mu}(\eta ) a(\eta ) b^{-s}(\eta ) \|_{0,0}
       \| b^s(\eta ) \|_{s,0}. 
\end{multline*}
Applying \eqref{opvasy} and Definition \ref{ordred} (iii)
we obtain \eqref{symbest} with $A = B(- s+\mu,-s+ \nu,0) + B(s,s,0)$, together with the uniform boundedness of the involved constants.
\qed \\
Also here it can be proved that the involved constants in Propositions \ref{732.a2201.pr}, \ref{s,s-nu} are uniform in compact sets with respect to $s,\mu,\nu$. 
\begin{Prop} \label{732.2201.pr}
The symbol spaces have the following properties:
\begin{enumerate}
\item[\textup{(i)}] $S^\mu(U\times\R^q;\E,\widetilde\E)\subseteq S^{\mu'}(U\times\R^q;\E,
\widetilde\E)$ for every $\mu'\geq\mu$;
\item[\textup{(ii)}] $D_y^\alpha D_\eta^\beta S^\mu(U\times\R^q;\E,\widetilde\E)\subseteq 
S^{\mu-|\beta|}(U\times\R^q;\E,\widetilde\E)$ for every $\alpha\in\N^p,\beta\in\N^q$;
\item[\textup{(iii)}] $S^\mu(U\times\R^q;\E_0,\widetilde\E)S^\nu(U\times\R^q;\E,\E_0)\subseteq
S^{\mu+\nu}(U\times\R^q;\E,\widetilde\E)$ for every $\mu,\nu\in\R$ \lr the notation on the left 
hand side of the latter relation means the space of all $(y,\eta)$-wise compositions of 
elements in the respective factors\rr.
\end{enumerate}
\end{Prop}
\proof
For simplicity we consider symbols with constant coefficients. 
Let us write $\|\cdot\|:=\|\cdot\|_{0,0}$, etc.

(i) $a(\eta)\in S^\mu(\R^q;\E,\widetilde\E)$ means \eqref{732.est1901.eq} and \eqref{opvasy};
      this implies
\begin{multline*}
  \| \td b^{s-\mu' +|\beta|}(\eta )  \{ D_\eta^\beta
        a(\eta ) \} b^{-s}(\eta )\|  =  \ \| \td b^{\mu-\mu'}(\eta ) \td b^{s- \mu +
       |\beta |}(\eta ) \{ D_\eta^\beta a(\eta ) \}
       b^{-s}(\eta ) \|  \\
   \leq  \ c \la\eta\ra^{\mu-\mu'} \| 
            \td b^{s-\mu + |\beta |}(\eta ) \{
	  D_\eta^\beta a(\eta ) \} b^{-s}(\eta ) \|
   \leq  \ c \| \td b^{s- \mu+ | \beta|}(\eta ) \{
         D_\eta^\beta a(\eta ) \} b^{-s}(\eta ) \|.	      	
\end{multline*}
      We employed $\mu - \mu' \leq 0$ and the property
      (iv) in Definition \ref{ordred}.	

(ii) The estimates \eqref{732.est1901.eq} can be 
      written as
      \[  \| \tilde{b}^{s-(\mu- |\beta |)}(\eta ) \{
            D_\eta^\beta a(\eta ) \} b^{-s}(\eta ) \| \leq c \]
      which just means that $D_\eta^\beta a(\eta) \in
      S^{\mu-|\beta|}(\R^q;\E,\widetilde\E)$.

(iii) Given $a(\eta) \in S^\mu(\R^q;\E_0,\widetilde\E), \tilde{a}(\eta) \in S^\nu(\R^q;
      \E, \E_0)$ we have (with obvious meaning of notation)
    \[  \| \tilde{b}_0^{s-\nu+|\gamma|}(\eta)  \{
         D_\eta^\gamma \tilde{a}(\eta ) \} b^{-s}
	 (\eta) \|  \leq c, \          	    
    \| \tilde{b}^{s-\mu+|\delta|}(\eta )  
         \{ D_\eta^\delta a(\eta ) \}
        b_0^{-s}(\eta) \| \leq c    \]	
     for all $\gamma, \delta \in \N^q$. If $\alpha \in
     \N^q$ is any multi-index, $D_\eta^\alpha(a
     \tilde{a})(\eta)$ is a linear combination of
     compositions $D_\eta^\delta a(\eta) D_\eta^\gamma
     \tilde{a}(\eta)$ with $| \gamma | + | \delta| = |
     \alpha |$. It follows that
    \begin{multline}  \label{new11}
   \|   \tilde{b}^{s-(\mu +\nu ) + | \alpha |} (\eta )
          D_\eta^\delta a(\eta ) \{ D_\eta^\gamma
	  \tilde{a}(\eta ) \} b^{-s}(\eta ) \|  \\
   =  \ \| \tilde{b}^{s-(\mu +\nu )+| \alpha |}(\eta )
           D_\eta^\delta a(\eta ) b_0^{-s+\nu - |\gamma |}(\eta )
           b_0^{s-\nu +|\gamma |} (\eta ) D_\eta^\gamma
           \tilde{a}(\eta ) b^{-s}(\eta ) \|  \\  	   	   
     \leq  \ \| \tilde{b}^{t-\mu +|\alpha |- | \gamma |}(\eta )
            D_\eta^\delta a(\eta ) b_0^{-t}(\eta ) \| \ 
	    \| b_0^{s-\nu +|\gamma |}(\eta) D_\eta^\gamma
	    \tilde{a}(\eta ) b^{-s}(\eta ) \| 
     \end{multline}
for $t = s- \nu + | \gamma |$; the right hand side is 
bounded in $\eta$, since $|\alpha | -|\gamma | 
= |\delta|$.
\qed
\begin{Rem}
Observe from \eqref{new11} that the semi-norms of compositions of symbols can be estimated by products of semi-norms of the factors.
\end{Rem}

\subsection{An example from the parameter-dependent cone calculus}
We now construct a specific family of reductions of orders between weighted spaces on a compact 
manifold $M$ with conical singularity $v$, locally near $v$ modelled on a cone 
\[ X^\Delta :=(\Ro_+\times X)/(\{0\}\times X) \]
with a smooth compact manifold $X$ as base.
The parameter $\eta$ will play the role of covariables of the calculus of operators on a manifold
with edge; that is why we talk about an example from the edge calculus. The associated ``abstract'' 
cone calculus according to what we did so far in the Sections \ref{OrdRed} and \ref{opvalsymb} and
then below in Chapter \ref{Oprefex} will be a contribution to the calculus of corner operators of second 
generation.
It will be convenient to pass to the stretched manifold $\M$ associated with $M$ which is a compact 
$C^\infty$ manifold with boundary $\partial\M\cong X$ such that when we squeeze down $\partial\M$ to a single point $v$ we just recover $M$. Close to $\partial\M$
the manifold $\M$ is equal to a cylinder $[0,1)\times X\ni(t,x)$, a collar 
neighbourhood of $\partial\M$ in $M$. A part of the considerations will be performed on the open stretched cone $X^\wedge:=\R_+\times X\ni (t,x)$ where we identify $(0,1)\times X$ with the interior of the collar neighbourhood (for convenience, without indicating any pull backs of functions or operators with respect to that identification). Let $\widetilde{M}:=2\M$ be the double of $\M$ (obtained by gluing together two copies $\M_\pm$ of $\M$ along the common boundary $\partial\M$, where we identify $\M$ with $\M_+$); then $\widetilde M$ is a closed compact $C^\infty$ manifold. On the space $M$ we have a family of weighted Sobolev spaces $H^{s,\gamma}(M), s,\gamma\in\R$, that may be defined as
\[ H^{s,\gamma}(M):=\{\sigma u+(1-\sigma)v:u\in\Hi^{s,\gamma}(X^\wedge),v\in H^s_\loc(M\setminus\{v\})\}, \]
where $\sigma(t)$ is a cut-off function (i.e., $\sigma\in C^\infty_0(\Ro_+), \sigma\equiv 1$ near $t=0$), $\sigma(t)=0$ for $t>2/3$. Here $\Hi^{s,\gamma}(X^\wedge)$ is defined to be the completion of $C_0^\infty(X^\wedge)$ with respect to the norm
\begin{equation} \label{neu14m}
 \left\{\frac{1}{2\pi i}\int_{\Gamma_{\frac{n+1}{2}-\gamma}}\|b^\mu_
\textup{base}(\imb\, w)(\Mu u)(w)\|^2_{L^2(X)}dw\right\}^\oneh ,
\end{equation}
$n=\dim X$, where $b_\textup{base}^\mu (\tau)\in L^\mu_\clw(X;\R_\tau )$ is a family of reductions of order on $X$, similarly as in Example \ref{ExamX} (in particular, $b^s_\textup{base}(\tau):H^s(X)\to H^0(X)=L^2(X)$ is a family  of isomorphisms). Moreover, $\Mu$ is the Mellin transform, $(\Mu u)(w)=\int_0^\infty t^{w-1}u(t)dt$, $w\in\C$ the complex Mellin covariable, and
\[ \Gamma_\beta:=\{w\in\C:\reb\, w=\beta\} \] 
for any real $\beta$. From $t^\delta\Hi^{s,\gamma}(X^\wedge)=\Hi^{s,\gamma+\delta}(X^\wedge)$ for all $s,\gamma,\delta\in\R$ it follows the existence of a strictly positive function $\textup{h}^\delta\in C^\infty(M\setminus\{v\})$, such that the operator of multiplication by $\textup{h}^\delta$ induces an isomorphism
\begin{equation} \label{equiso}
\textup{h}^\delta:H^{s,\gamma}(M)\to H^{s,\gamma+\delta}(M)
\end{equation}
for every $s,\gamma,\delta\in\R$.\\
Moreover, again according to Example \ref{ExamX}, now for any smooth compact manifold $\widetilde M$ we have an order reducing family $\td b(\eta)$ in the scale of Sobolev spaces $H^s(\widetilde M),s\in\R$. More generally, we employ parameter-dependent families $\td a(\eta)\in L_\clw^\mu(\widetilde M;\R^q)$. The symbols $a(\eta)$ that we want to establish in the scale $H^{s,\gamma}(M)$ on our compact manifold $M$ with conical singularity $v$ will be essentially (i.e., modulo Schwartz functions in $\eta$ with values in globally smoothing operators on $M$) constructed in the form
\begin{equation} \label{bmu}
a(\eta):=\sigma a_\textup{edge}(\eta)\td\sigma +(1-\sigma) a_\textup{int}(\eta) (1-\tilde{\tilde{\sigma}}),
\end{equation} 
$a_\textup{int} (\eta):=\td a(\eta)|_{\textup{int}\M}$, with cut-off functions $\sigma(t),\td\sigma(t),\tilde{\tilde{\sigma}}(t)$ on the half axis, supported in $[0,2/3)$, with the property
\[ \tilde{\tilde{\sigma}}\prec\sigma\prec\td\sigma \]
(here $\sigma\prec\td\sigma$ means the $\td\sigma$ is equal to $1$ in a neighbourhood of supp\, $\sigma$).\\ 
The ``edge'' part of \eqref{bmu} will be defined in the variables $(t,x)\in X^\wedge$. Let us choose a parameter-dependent elliptic family of operators of order $\mu$ on $X$
\[ \td p(t,\td\tau,\td\eta)\in C^\infty(\Rb_+,L^\mu_\clw(X;\R^{1+q}_{\td\tau,\td\eta})). \]
Setting 
\begin{equation} \label{rhoeq}
p(t,\tau,\eta):=\td p(t,t\tau,t\eta)
\end{equation}
we have what is known as an edge-degenerate family of operators on $X$. We now employ the following Mellin quantisation theorem. 
\begin{Def}
Let $M_\mathcal{O}^\mu(X;\R^q)$ defined as the set of all $h(z,\eta)\in \mathcal{A} (\C, L_\clw^\mu(X;\R^q))$ such that $h(\beta+i\tau,\eta)\in L_\clw^\mu(X;\R_{\tau,\eta}^{1+q})$ for every $\beta\in\R$, uniformly in compact $\beta$-intervals 
\lr here $\mathcal{A}(\C,E)$ with any Fr\'echet space $E$ denotes the space of all $E$-valued holomorphic functions in $\C$, in the Fr\'echet topology of uniform convergence on compact sets\rr.
\end{Def}
Observe that also $M_\mathcal{O}^\mu(X;\R^q)$ is a Fr\'echet space in a natural way. Given an $f(t,t',z,\eta)\in C^\infty(\R_+\times\R_+,L^\mu_\clw(X;\Gamma_{\oneh-\gamma}\times\R^q))$ we set
\[ \sop_M^\gamma(f)(\eta)u(r):=\int_\R\int_0^\infty (\frac{t}{t'})^{-(\oneh-\gamma+i\tau)}f(t,t',\oneh-\gamma+i\tau,\eta)u(t')\frac{dt'}{t'}\dbar\tau \]
which is regarded as a (parameter-dependent) weighted pseudo-differential operator with symbol $f$, referring to the weight 
$\gamma\in\R$. There exists an element
\begin{equation} \label{h1eq}
\td h(t,z,\td\eta)\in C^\infty(\Rb_+,M^\mu_\mathcal{O} (X;\R^{q}_{\td\eta})) 
\end{equation}
such that, when we set 
\begin{equation} \label{heq}
h(t,z,\eta):=\td h(t,z,t\eta) 
\end{equation}
we have
\begin{equation} \label{newe15}
\sop_M^\gamma (h)(\eta)=\op_t(p)(\eta) 
\end{equation}
mod $L^{-\infty}(X^\wedge;\R^{q}_{\eta})$, for every weight $\gamma\in\R$. Observe that when we set 
\[ p_0(t,\tau,\eta):=\td p(0,t\tau,t\eta), \ h_0(t,z,\eta):=\td h(0,z,t\eta) \]
we also have $\sop_M^\gamma(h_0)(\eta)=\op_t(p_0)(\eta) \mod L^{-\infty}(X^\wedge;\Rb_+)$, for all $\gamma\in\R$.\\
Let us now choose cut-off functions $\omega(t),\td\omega(t),\tilde{\tilde{\omega}}(t)$ such that $\tilde{\tilde{\omega}}\prec\omega\prec\td\omega$.\\
Fix the notation $\omega_\eta(t):=\omega(t[\eta])$, and form the operator function
\begin{multline} \label{bedgeex}
a_\textup{edge}(\eta):= \omega_\eta (t)t^{-\mu} \sop_M^{\gamma-\frac{n}{2}}(h)(\eta)\td\omega_\eta(t) \\  +t^{-\mu}\big(1-\omega_\eta(t)\big)\op_t(p)
(\eta)\big(1-\tilde{\tilde{\omega}}_\eta(t)\big) 
+m(\eta)+g(\eta). 
\end{multline}
Here $m(\eta)$ and $g(\eta)$ are smoothing Mellin and Green symbols of the edge calculus.
The definition of $m(\eta)$ is based on smoothing Mellin symbols $f(z)\in M^{-\infty}(X;\Gamma_\beta)$. Here $M^{-\infty}(X;\Gamma_\beta)$ is the subspace of all $f(z)\in L^{-\infty}(X;\Gamma_\beta)$ such that for some $\varepsilon >0$ (depending on $f$) the function $f$ extends to an 
\[ l(z)\in\An(U_{\beta,\varepsilon},L^{-\infty}(X)) \]
where $U_{\beta,\varepsilon}:=\{z\in\C:|\reb z-\beta|<\varepsilon\}$ and
\[ l(\delta+i\tau)\in L^{-\infty}(X;\R_\tau) \]
for every $\delta\in(\beta-\varepsilon,\beta+\varepsilon)$, uniformly in compact subintervals. By definition we then have $f(\beta+i\tau)=l(\beta+i\tau)$; for brevity we often denote the holomorphic extension $l$ of $f$ again by $f$. For $f\in M^{-\infty}(X;\Gamma_{\frac{n+1}{2}-\gamma})$ we set 
\[ m(\eta):=t^{-\mu}\omega_\eta\sop_M^{\gamma-\frac{n}{2}}(f)\td\omega_\eta \]
for any cut-off  functions $\omega,\td\omega$.\\ 
In order to explain the structure of $g(\eta)$ in \eqref{bedgeex} we first introduce weighted spaces on the infinite stretched cone $X^\wedge=\R_+\times X$, namely,
\begin{equation} \label{jadid16} \K^{s,\gamma;g}(X^\wedge):=\omega\Hi^{s,\gamma}(X^\wedge)+(1-\omega)H^{s;g}_\cone(X^\wedge) 
\end{equation}
for any $s,\gamma,g\in\R$, and a cut-off function $\omega$, see \eqref{neu14m} which defines $\Hi^{s,\gamma}(X^\wedge)$ and the formula \eqref{equtcon}. Moreover, we set $\K^{s,\gamma}(X^\wedge):=\K^{s,\gamma;0}(X^\wedge)$. The operator families $g(\eta)$ are so-called Green symbols in the covariable $\eta\in\R^q$, defined by 
\begin{equation} \label{equt161} 
g(\eta)\in S_\clw^\mu(\R^q_\eta;\K^{s,\gamma;g}(X^\wedge),\Sch^{\gamma-\mu+\varepsilon} (X^\wedge)),
\end{equation}
\begin{equation} \label{equt162}
g^*(\eta)\in S_\clw^\mu(\R^q_\eta;\K^{s,-\gamma+\mu;g}(X^\wedge),\Sch^{-\gamma+\varepsilon} (X^\wedge)),
\end{equation}
for all $s,\gamma,g\in\R$, where $g^*$ denotes the $\eta-$wise formal adjoint with respect to the scalar product of $\K^{0,0;0}(X^\wedge)=r^{-\frac{n}{2}}L^2(\R_+\times X)$ and $\varepsilon =\varepsilon(g)>0$. Here 
\[ \Sch^\beta(X^\wedge):=\omega \K^{\infty,\beta}(X^\wedge)+(1-\omega)\Sch(\Rb_+,C^\infty(X)) \]
for any cut-off function $\omega$. The notion of operator-valued symbols in \eqref{equt161}, \eqref{equt162} refers to \eqref{neu5} in its generalisation to Fr\'echet spaces $\widetilde H$ (rather than Hilbert spaces) with group actions (see Remark \ref{rem7}) that is in the present case given by 
\begin{equation} \label{eut+}
\kappa_\lambda: u(t,x)\to\lambda^{\frac{n+1}{2}+g}u(\lambda t,x),\ \lambda\in\R_+
\end{equation}
$n=\dim X$, both in the spaces $\K^{s;\gamma,g}(X^\wedge)$ and $\Sch^{\gamma-\mu+\varepsilon}(X ^\wedge)$.\\
The following Theorem \ref{thmsymb} is crucial for proving that our new order reduction family is well defined. Therefore we will sketch the main steps of the proof, which is based on the edge calculus. Various aspects of the proof can be found in the literature, for example in Kapanaze and Schulze \cite[Proposition 3.3.79]{Kapa5}, Schrohe and Schulze \cite{Schr19}, Harutyunyan and Schulze \cite{Haru13}. Among the tools we have the pseudo-differential operators on $X^\wedge$ interpreted as a manifold with conical exit to infinity $r\to\infty$; the general background may be found in Schulze \cite{Schu20}. The calculus of such exit operators goes back to Parenti \cite{Pare1}, Cordes \cite{Cord1}, Shubin \cite{Shub1}, and others.\\
\begin{Thm} \label{thmsymb}
We have 
\begin{equation} \label{newe16} 
\sigma a_\textup{edge}(\eta)\td\sigma\in S^\mu(\R^q;\K^{s,\gamma;g}(X^\wedge),\K^{s-\mu,\gamma -\mu;g}(X^\wedge))
\end{equation} 
for every $s,g\in\R$, more precisely,
\begin{equation} \label{equ+}
D_\eta^\beta \{\sigma a_\textup{edge}(\eta)\td\sigma\}\in S^{\mu-|\beta|}(\R^q;\K^{s,\gamma;g} (X^\wedge),\K^{s-\mu+|\beta|,\gamma-\mu;g}(X^\wedge)) 
\end{equation}
for all $s,g\in\R$ and all $\beta\in\N^q$. \lr The spaces of symbols in \eqref{newe16}, \eqref{equ+} refer to the group action \eqref{eut+}\rr.
\end{Thm}
\proof
To prove the assertions it is enough to consider the case without $m(\eta)+g(\eta)$, since the latter sum maps to $\K^{\infty,\gamma;g}(X^\wedge)$ anyway.
The first part of the Theorem is known, see, for instance, \cite{Haru13} or \cite{Dine5}.
Concerning the relation \eqref{equ+} we write
\begin{equation} \label{equtsq}
\sigma a_\textup{edge}(\eta)\td\sigma =\sigma\{ a_\textup{c}(\eta)+ a_\psi(\eta)\}\td\sigma 
\end{equation}
with
\[ a_\textup{c}(\eta):=t^{-\mu}\omega_\eta\sop_M^{\gamma-\frac{n}{2}}(h)(\eta)\td\omega_\eta , \]
\[ a_\psi(\eta):=t^{-\mu}(1-\omega_\eta)\op_t(p)(\eta)(1-\td{\td\omega}_\eta)\]
and it suffices to take the summands separately.
In order to show \eqref{equ+} we consider, for instance, the derivative $\partial /\partial\eta_j =:\partial_j$ for some $1\leq j\leq q$. By iterating the process we then obtain the assertion. We have 
\[\partial_j \sigma\{a_\textup{c}(\eta)+a_\psi(\eta)\}\td\sigma=\sigma\{\partial_j a_\textup{c}(\eta)+\partial_j a_\psi (\eta)\}\td\sigma=b_1(\eta)+b_2(\eta)+b_3(\eta) \]
with
\begin{align*}
b_1(\eta) &:=\sigma t^{-\mu}\Big\{\omega_\eta\sop_M^{\gamma-\frac{n}{2}}(h)(\eta)\partial_j \td\omega_\eta + (1-\omega_\eta)\op_t(p)(\eta)\partial_j(1-\td{\td\omega}_\eta)\Big\}\td\sigma, \\
b_2(\eta) &:=\sigma t^{-\mu}\Big\{\omega_\eta \sop_M^{\gamma-\frac{n}{2}}(\partial_j h)(\eta)\td \omega_\eta +(1-\omega_\eta)\op_t(\partial_j p)(\eta)(1-\td{\td\omega}_\eta)\Big\} \td{\td\sigma}, \\
b_3(\eta) &:=\sigma t^{-\mu}\Big\{(\partial_j\omega_\eta)\sop_M^{\gamma-\frac{n}{2}}(h)(\eta) \td\omega_\eta + (\partial_j(1-\omega_\eta))\op_t(p)(\eta)(1-\td{\td\omega}_\eta)\Big\}\td\sigma .
\end{align*}
In $b_1(\eta)$ we can apply a pseudo-locality argument which is possible since $\partial_j \td\omega_\eta\equiv 0$ on $\textup{supp}\,\omega_\eta$ and $\partial_j(1-\td{\td\omega}_\eta)\equiv 0$ on $\textup{supp}\,(1-\omega_\eta)$; this yields (together with similar considerations as for the proof of \eqref{newe16})
\[b_1(\eta)\in S^{\mu -1}(\R^q;\K^{s,\gamma;g}(X^\wedge),\K^{\infty,\gamma-\mu;g}(X^\wedge)). \]
Moreover we obtain
\[ b_2(\eta)\in S^{\mu-1}(\R^q;\K^{s,\gamma;g}(X^\wedge),\K^{s-\mu+1,\gamma-\mu;g}(X^\wedge)) \]
since $\partial_j h$ and $\partial_j p$ are of order $\mu-1$ (again combined with arguments for \eqref{newe16}). Concerning $b_3(\eta)$ we use the fact that there is a $\psi\in C_0^\infty (\R_+)$ such that $\psi\equiv 1$ on $\textup{supp}\,\partial_j\omega$, $\td\omega-\psi\equiv 0$ on $\textup{supp}\,\partial_j\omega$ and $(1-\td{\td\omega})-\psi\equiv 0$ on $\textup{supp}\, \partial_j\omega$. Thus, when we set $\psi_\eta(t):=\psi(t[\eta])$, we obtain $b_3(\eta):=c_3(\eta)+c_4(\eta)$ with
\begin{align*}
c_3(\eta) &:=\sigma t^{-\mu}\Big\{(\partial_j\omega_\eta)\sop_M^{\gamma-\frac{n}{2}}(h)(\eta)\psi_ \eta-(\partial_j\omega_\eta)\op_t(p)(\eta)\psi_\eta\Big\}\td\sigma, \\
c_4(\eta) &:=\sigma t^{-\mu}\Big\{(\partial_j\omega_\eta)\sop_M^{\gamma-\frac{n}{2}}(h)(\eta)[\td \omega_\eta-\psi_\eta] -(\partial_j\omega_\eta)\op_t(p)(\eta) [(1-\td{\td\omega}_\eta)-\psi_\eta] \Big\}\td\sigma.
\end{align*}
Here, using $\partial_j\omega_\eta=(\omega')_\eta\partial_j(t[\eta])$ which yields an extra power of $t$ on the left of the operator, together with pseudo-locality, we obtain
\[ c_4(\eta)\in S^{\mu-1}(\R^q;\K^{s,\gamma;g}(X^\wedge),\K^{\infty,\gamma-\mu;g}(X^\wedge)). \]
To treat $c_3(\eta)$ we employ that both $\partial_j\omega_\eta$ and $\psi_\eta$ are compactly supported on $\R_+$. Using the property \eqref{newe15}, we have 
\begin{multline*} 
c_3(\eta)=\sigma t^{-\mu}(\partial_j\omega_\eta)\big\{\sop_M^{\gamma-\frac{n}{2}}(h)(\eta)-\op_t (p)(\eta)\big\}\psi_\eta\td\sigma \\ 
\in S^{\mu-1}(\R^q;\K^{s,\gamma;g}(X^\wedge),\K^{\infty,\gamma-\mu; g}(X^\wedge)).
\end{multline*}
\qed
\begin{Def}
A family of operators $c(\eta)\in\Sch(\R^q,\bigcap_{s\in\R}\Li(H^{s,\gamma}(M),H^{\infty,\delta} (M)))$ is called a smoothing element in the parameter-dependent cone calculus on $M$ associated with the weight data $(\gamma,\delta)\in\R^2$, written $c\in C_G(M,(\gamma,\delta);\R^q)$, if there is an $\varepsilon=\varepsilon(c) >0$ such that
\[ \begin{array}{c} 
c(\eta)\in \Sch(\R^q,\Li(H^{s,\gamma}(M),H^{\infty,\delta+\varepsilon}(M))), \\ 
c^*(\eta)\in \Sch(\R^q,\Li(H^{s,-\delta}(M),H^{\infty,-\gamma+\varepsilon}(M))); 
\end{array}\]
for all $s\in\R$; here $c^*$ is the $\eta$-wise formal adjoint of $c$ with respect to the $H^{0,0}(M)$-scalar product.
\end{Def}
The $\eta$-wise kernels of the operators $c(\eta)$ are in $C^\infty\left( (M\setminus\{v\})\times(M\setminus\{v\})\right)$. However, they are of flatness $\varepsilon$ in the respective distance variables to $v$, relative to the weights $\delta$ and $\gamma$, respectively. Let us look at a simple example to illustrate the structure. We choose elements $k\in \Sch(\R^q, H^{\infty,\delta +\varepsilon}(M))$, $k'\in \Sch(\R^q,H^{\infty,-\gamma+\varepsilon}(M))$ and assume for convenience that $k$ and $k'$ vanish outside a neighbourhood of $v$, for all $\eta\in\R^q$. Then with respect to a local splitting of variables $(t,x)$ near $v$ we can write $k=k(\eta,t,x)$ and $k'=k'(\eta,t',x')$, respectively. Set 
\[ c(\eta)u(t,x):=\iint k(\eta,t,x)k'(\eta,t',x')u(t',x')t'^n dt'dx' \]
with the formal adjoint
\[ c^*(\eta)v(t',x'):=\iint \overline{k'(\eta,t',x')k(\eta,t,x)}v(t,x)t^n dtdx. \]
Then $c(\eta)$ is a smoothing element in the parameter-dependent cone calculus.\\
By $C^\mu(M,(\gamma,\gamma-\mu);\R^q)$ we denote the set of all operator families 
\begin{equation} \label{qe23}
a(\eta)=\sigma a_\textup{edge}(\eta)\td\sigma+(1-\sigma)a_\textup{int}(\eta)(1-\td{\td\sigma})+ c(\eta)
\end{equation}
where $a_\textup{edge}$ is of the form \eqref{bedgeex}, $a_\textup{int}\in L_\clw^\mu(M\setminus\{v\} ;\R^q)$, while $c(\eta)$ is a parameter-dependent smoothing operator on $M$, associated with the weight data $(\gamma,\gamma-\mu)$.
\begin{Thm} \label{thm132}
Let $M$ be a compact manifold with conical singularity. Then the $\eta$-dependent families \eqref{bmu} which define continuous operators 
\begin{equation} \label{conaedge}
a(\eta):H^{s,\gamma}(M)\to H^{s-\nu,\gamma-\nu}(M) 
\end{equation}
for all $s\in\R$, $\nu\geq\mu$, have the properties:
\begin{equation} \label{equatB}
\| a(\eta)\|_{\Li(H^{s,\gamma}(M),H^{s-\nu,\gamma-\nu}(M))}\leq c\la\eta\ra^B 
\end{equation}
for all $\eta\in\R^q$, and $s\in\R$, with constants $c=c(\mu,\nu,s)>0$, $B=B(\mu,\nu,s)$, and, when $\mu\leq 0$
\begin{equation} \label{equatm}
\| a(\eta) \|_{\Li(H^{0,0}(M),H^{0,0}(M))}\leq c\la\eta\ra^\mu 
\end{equation}
for all $\eta\in\R$, $s\in\R$, with constants $c=c(\mu,s)>0$. 
\end{Thm}
\proof
The result is known for the summand $(1-\sigma)a_\textup{int}(\eta)(1-\td{\td\sigma})$ as we see from Example \ref{ExamX}. Therefore, we may concentrate on 
\[p(\eta):=\sigma a_\textup{edge}(\eta)\td\sigma:H^{s,\gamma}(M)\to H^{s-\nu,\gamma-\nu}(M).\]
To show \eqref{equatB} we pass to
\[\sigma a_\textup{edge}(\eta)\td\sigma:\K^{s,\gamma}(X^\wedge)\to\K^{s-\nu,\gamma-\nu}
(X^\wedge).\]
Then Theorem \ref{thmsymb} shows that we have symbolic estimates, especially
\[ \|\kappa^{-1}_{\la\eta\ra}p(\eta)\kappa_{\la\eta\ra}\|_{\Li(\K^{s,\gamma}(X^\wedge), \K^{s-\mu,\gamma-\mu}(X^\wedge))}\leq c\la\eta\ra^\mu. \]
We have
\[ \| p(\eta)\|_{\Li(\K^{s,\gamma}(X^\wedge), \K^{s-\nu,\gamma-\nu}(X^\wedge))} \leq 
   \| p(\eta)\|_{\Li(\K^{s,\gamma}(X^\wedge), \K^{s-\mu,\gamma-\mu}(X^\wedge))}, \]
 and
\begin{multline} \nonumber
\| p(\eta)\|_{\Li(\K^{s,\gamma}(X^\wedge), \K^{s-\mu,\gamma-\mu}(X^\wedge))}=
\| \kappa_{\la\eta\ra}\kappa^{-1}_{\la\eta\ra}p(\eta)\kappa_{\la\eta\ra}\kappa^{-1}_{\la\eta\ra} \|_
{\Li(\K^{s,\gamma}(X^\wedge), \K^{s-\mu,\gamma-\mu}(X^\wedge))} \\
\leq \|\kappa_{\la\eta\ra}\|_{\Li(\K^{s-\mu,\gamma-\mu}(X^\wedge), \K^{s-\mu,\gamma-\mu}(X^\wedge))} \|\kappa^{-1}_{\la\eta\ra}p(\eta)\kappa_{\la\eta\ra}\|_{\Li(\K^{s,\gamma}(X^\wedge), \K^{s,\gamma}(X^\wedge))} \\ \|\kappa^{-1}_{\la\eta\ra}\|_{\Li(\K^{s-\mu,\gamma-\mu}(X^\wedge), \K^{s,\gamma}(X^\wedge))} \leq 
c\la\eta\ra^{\mu+\widetilde M+M}.
\end{multline}
Here we used that $\kappa_{\la\eta\ra},\kappa_{\la\eta\ra}^{-1}$ satisfy estimates like \eqref{kappMc}. \\
For \eqref{equatm} we employ that $\kappa_\lambda$ is operating as a unitary group on $\K^{0,0}(X^\wedge)$. This gives us 
\begin{multline} \nonumber
\|p(\eta)\|_{\Li(\K^{0,0}(X^\wedge),\K^{0,0}(X^\wedge))}=\|\kappa^{-1}_{\la\eta\ra}p(\eta)\kappa_
{\la\eta\ra}\|_{\Li(\K^{0,0}(X^\wedge),\K^{0,0}(X^\wedge))} \\
\leq \|\kappa^{-1}_{\la\eta\ra}p(\eta)\kappa_
{\la\eta\ra}\|_{\Li(\K^{0,0}(X^\wedge),\K^{-\mu,-\mu}(X^\wedge))}\leq c\la\eta\ra^\mu.
\end{multline}
\qed
\begin{Thm} \label{thmind15}
For every $k\in\Z$ there exists an $f_k(z)\in M^{-\infty}(X;\Gamma_{\frac{n+1}{2}-\gamma})$ such that for every cut-off functions 
$\omega,\td\omega$ the operator
\begin{equation} \label{neww23} 
A:=1+\omega\sop_M^{\gamma-\frac{n}{2}}(f_k)\td\omega:H^{s,\gamma}(M)\to H^{s,\gamma}(M) 
\end{equation}
is Fredholm and of index $k$, for all $s\in\R$.
\end{Thm}
\proof
We employ the result (cf. \cite{Schu31}) that for every $k\in\Z$ there exists an $f_k(z)$ such that 
\begin{equation} \label{nnew23}
\widetilde A:=1+\omega\sop_M^{\gamma-\frac{n}{2}}(f_k)\td\omega:\K^{s,\gamma}(X^\wedge)\to \K^{s,\gamma}(X^\wedge)
\end{equation}
is Fredholm of index $k$. Recall that the proof of the latter result follows from a corresponding theorem in the case $\dim X=0$. The Mellin symbol $f_k$ is constructed in such a way that $1+f_k(z)\not=0$ for all $z\in\Gamma_{\oneh-\gamma}$ and the argument of $1+f_k(z)|_{\Gamma_{\oneh-\gamma}}$ varies from 1 to $2\pi k$ when $z\in\Gamma_{\oneh-\gamma}$ goes from $\imb z=-\infty$ to $\imb z=+\infty$. The choice of $\omega,\td\omega$ is unessential; so we assume that $\omega,\td\omega\equiv 0$ for $r\geq 1-\varepsilon$ with some $\varepsilon >0$. Let us represent the cone $\widetilde M:=X^\Delta$ as a union of $\big([0,1+\frac{\varepsilon}{2}) \times X\big)/(\{0\} \times X)=:\widetilde M_-$ and $(1-\frac{\varepsilon}{2},\infty)\times X=:\widetilde M_+$. Then 
\begin{equation} \label{tildeA}
\widetilde A|_{\widetilde M_-}=1+\omega\sop_M^{\gamma-\frac{n}{2}}(f_k)\td\omega, \ \widetilde A|_{\widetilde M_+}=1.
\end{equation}
Moreover, without loss of generality, we represent $M$ as a union $\big([0,1+\frac{\varepsilon} {2}) \times X\big)/ (\{0\}\times X)\cup M_+$ where $M_+$ is an open $C^\infty$ manifold which intersects $\big([0,1+\frac{\varepsilon} {2}) \times X\big)/ (\{0\}\times X)=:M_-$ in a cylinder of the form $(1-\frac{\varepsilon}{2},1+\frac{\varepsilon}{2})\times X$.
Let $B$ denote the operator on $M$, defined by 
\begin{equation} \label{opA} 
B_-:=A|_{M_-}=1+\omega\sop_M^{\gamma-\frac{n}{2}}(f_k)\td\omega, \ B_+:=A|_{M_+}=1
\end{equation}
We are then in a special situation of cutting and pasting of Fredholm operators. We can pass to manifolds with conical singularities $N$ and $\widetilde N$ by setting 
\[ N=\widetilde M_-\cup M_+, \ \widetilde N=M_-\cup \widetilde M_+ \]
and transferring the former operators in \eqref{tildeA}, \eqref{opA} to $N$ and $\widetilde N$, respectively, by gluing together the $\pm$ pieces of $\widetilde A$ and $A$ to belong to $\widetilde M_\pm$ and $M_\pm$ to corresponding operators $\widetilde B$ on $\widetilde N$ and $B$ on $N$. We then have the relative index formula
\begin{equation} \label{indalph}
\ind A-\ind B=\ind \widetilde A-\ind \widetilde B
\end{equation}
(see \cite{Naza10}). In the present case $\widetilde A$ and $\widetilde M$ are the same as $B$ and $N$ where $\widetilde B$ and $\widetilde N$ are the same as $A$ and $M$. It follows that 
\begin{equation}\label{indbet}
\ind \widetilde A-\ind\widetilde B=\ind B-\ind A.
\end{equation}
From \eqref{indalph}, \eqref{indbet} it follows that $\ind A=\ind B=\ind\widetilde A$.
\qed
\begin{Thm} \label{thmti}
There is a choice of $m$ and $g$ such that the operators \eqref{bmu} form a family of isomorphisms 
\begin{equation}  \label{nnew24}
a(\eta):H^{s,\gamma}(M)\to H^{s-\mu,\gamma-\mu}(M) 
\end{equation}
for all $s\in\R$ and all $\eta\in\R^q$.
\end{Thm} 
\proof
We choose a function 
\[ p(t,\tau,\eta,\zeta):=\td p(t\tau,t\eta,\zeta) \]
similarly as \eqref{rhoeq} where $\td p(\td\tau,\td\eta,\zeta)\in L_\clw^\mu(X;\R^{1+q+l}_{\td\tau, \td\eta,\zeta}),l\geq 1,$ is a parameter-dependent elliptic with parameters $\td\tau,\td\eta,\zeta$. For purposes below we specify $\td p(t,\td\tau,\td\eta,\zeta)$ in such a way that the parameter-dependent homogeneous principal symbol in $(t,x,\td\tau,\xi,\td\eta,\zeta)$ for $(\td\tau,\xi,\td\eta,\zeta) \not=0$ is equal to
\[ (|\td\tau|^2+|\xi|^2+|\td\eta|^2+|\zeta|^2)^\frac{\mu}{2}.\] 
We now form an element 
\[ \td h(t,z,\td\eta,\zeta)\in M_\mathcal{O}^\mu(X;\R^{q+l}_{\td\eta,\zeta}) \]
analogously as \eqref{h1eq} such that
\[ h(t,z,\eta,\zeta):=\td h(t,z,t\eta,\zeta) \]
satisfies 
\[ \sop_M^\gamma(h)(\eta,\zeta)=\op_t(p)(\eta,\zeta) \]
$\mod L^{-\infty}(X^\wedge;\R^{q+l}_{\eta,\zeta})$. For every fixed $\zeta\in\R^l$ this is exactly as before, but in this way we obtain corresponding $\zeta$-dependent families of such objects. It follows 
\[ \sigma b_\textup{edge}(\eta,\zeta)\td\sigma=t^{-\mu}\sigma \left\{ \omega_\eta\sop_M^ {\gamma- \frac{n}{2}}(h)(\eta,\zeta)\td\omega_\eta+\chi_\eta\op_t(p)(\eta,\zeta)\td\chi_\eta \right\} \td\sigma \]
with
\[ \chi_\eta(t):=1-\omega_\eta(t), \ \td\chi_\eta(t):=1-\td{\td\omega}_\eta(t). \]
Let us form the principal edge symbol
\[ \sigma_\wedge(\sigma b_\textup{edge}\td\sigma)(\eta,\zeta)=t^{-\mu}\left\{ \omega_{|\eta|}\sop_M^{\gamma-\frac{n}{2}}(h)(\eta,\zeta)\td\omega_{|\eta|}+\chi_{|\eta|}\op_t(p) (\eta,\zeta)\td\chi_{|\eta|}\right\} \]
for $|\eta|\not=0$ which gives us a family of continuous operators
\begin{equation} \label{Feqt}
\sigma_\wedge(\sigma b_\textup{edge}\td\sigma)(\eta,\zeta): \K^{s,\gamma;g}(X^\wedge)
\to\K^{s-\mu,\gamma-\mu;g}(X^\wedge)
\end{equation}
which is elliptic as a family of classical pseudo-differential operators on $X^\wedge$. In addition it is exit elliptic on $X^\wedge$ with respect to the conical exit of $X^\wedge$ to infinity. In order that \eqref{Feqt} is Fredholm for the given weight $\gamma\in\R$ and all $s,g\in\R$ it is necessary and sufficient that the subordinate conormal symbol
\[ \sigma_\textup{c}\sigma_\wedge(\sigma b_\textup{edge}\td\sigma)(z,\zeta):H^s(X)\to H^{s-\mu}(X) \]
is a family of isomorphisms for all $z\in\Gamma_{\frac{n+1}{2}-\gamma}$. This is standard information from the calculus on the stretched cone $X^\wedge$. By definition the conormal symbol is just
\begin{equation}\label{cosyis}
\td h(0,z,0,\zeta):H^s(X)\to H^{s-\mu}(X).
\end{equation}
Since by construction $\td h(\beta+i\tau,0,\zeta)$ is parameter-dependent elliptic on $X$ with parameters $(\tau,\zeta)\in\R^{1+l}$, for every $\beta\in\R$ (uniformly in finite $\beta$-intervals) there is a $C>0$ such that \eqref{cosyis} becomes bijective whenever $|\tau,\zeta|>C$. In particular, choosing $\zeta$ large enough it follows the bijectivity for all $\tau\in\R$, i.e., for all $z\in\Gamma_{\frac{n+1}{2}-\gamma}$. Let us fix $\zeta^1$ in that way and write again
\[ p(t,\tau,\eta):=p(t,\tau,\eta,\zeta^1), \ h(t,z,\eta):=h(z,t\eta,\zeta^1). \]
We are now in the same situation we started with, but we know in addition that \eqref{Feqt} is a family of Fredholm operators of a certain index, say, $-k$ for some $k\in\Z$. With the smoothing Mellin symbol $f_k(z)$ as in \eqref{nnew23} we now form the composition
\begin{equation}
\label{Ceq}
\sigma b_\textup{edge}(\eta)\td\sigma(1+\omega_\eta\sop_M^{\gamma-\frac{n}{2}}(f_k)\td\omega_\eta)
\end{equation}
which is of the form
\begin{equation}
\label{Feq}
\sigma b_\textup{edge}(\eta)\td\sigma+\omega_\eta\sop_M^{\gamma-\frac{n}{2}}(f)\td\omega_\eta+g(\eta)
\end{equation}
for another smoothing Mellin symbol $f(z)$ and a certain Green symbol $g(\eta)$. Here, by a suitable choice of $\omega,\td\omega$, without loss of generality we assume that $\sigma\equiv 1$ and $\td\sigma\equiv 1$ on $\textup{supp}\, \omega_\eta\cup\textup{supp}\, \td\omega_\eta$, for all $\eta\in\R^q$.
Since \eqref{Ceq} is a composition of parameter-dependent cone operators the associated edge symbol is equal to 
\begin{equation}
\label{Feq}
F(\eta):=\sigma_\wedge(\sigma b_\textup{edge}\td\sigma)(\eta)(1+\omega_{|\eta|}\sop_M^{\gamma- \frac{n}{2}}(f_k)\td\omega_{|\eta|}):\K^{s,\gamma}(X^\wedge)\to \K^{s-\mu,\gamma-\mu}(X^\wedge)
\end{equation}
which is a family of Fredholm operators of index $0$. By construction \eqref{Feq} depends only on $|\eta|$. For $\eta\in S^{q-1}$ we now add a Green operator $g_0$ on $X^\wedge$ such that
\[ F(\eta)+g_0(\eta):\K^{s,\gamma}(X^\wedge)\to \K^{s-\mu,\gamma-\mu}(X^\wedge) \]
is an isomorphism; it is known that such $g_0$ (of finite rank) exists (for $N=\dim\ker F(\eta)$ it can be written in the form $g_0 u:=\sum_{j=1}^N(u,v_j)w_j$, where $(\cdot,\cdot)$ is the $K^{0,0}(X^\wedge)$-scalar product and $(v_j)_{j=1,\dots,N}$ and $(w_j)_{j=1,\dots,N}$ are orthonormal systems of functions in $C_0^\infty(X^\wedge)$). Setting
\[ g(\eta):=\sigma\vartheta(\eta)|\eta|^\mu\kappa_{|\eta|}g_0\kappa^{-1}_{|\eta|}\td\sigma \]
with an excision  function $\vartheta(\eta)$ in $\R^q$ we obtain a Green symbol with $\sigma_\wedge (g)(\eta)=|\eta|^\mu\kappa_{|\eta|}g_0\kappa^{-1}_{|\eta|}$ and hence
\[ \sigma_\wedge(F(\eta)+g(\eta)):\K^{s,\gamma}(X^\wedge)\to \K^{s-\mu,\gamma-\mu}(X^\wedge) \]
is a family of isomorphisms for all $\eta\in\R^q\setminus\{0\}$. Setting
\begin{multline} a_\textup{edge}(\eta):=\left[t^{-\mu}\omega_\eta\sop_M^{\gamma-\frac{n}{2}}(h)\td\omega_\eta+ \chi_\eta\op_t(p)(\eta)\td\chi_\eta\right]\left(1+\omega_\eta\sop_M^{\gamma-\frac{n}{2}}(f_k)\td \omega_\eta\right) \\ +|\eta|^\mu\vartheta(\eta)\kappa_{|\eta|}g_0\kappa^{-1}_{|\eta|} 
\end{multline}
we obtain an operator family
\[ \sigma a_\textup{edge}(\eta)\td\sigma=F(\eta)+g(\eta) \]
as announced before.
Next we choose a parameter-dependent elliptic $a_\textup{int}(\eta)\in L^\mu_\clw(M\setminus\{v\}; \R^q_\eta)$ such that its parameter-dependent homogeneous principal symbol close to $t=0$ (in the splitting of variables $(t,x)$) is equal to
\[ (|\tau|^2+|\xi|^2+|\eta|^2)^\frac{\mu}{2}. \]
 Then we form 
\[ a(\eta):=\sigma a_\textup{edge}(\eta)\td\sigma+(1-\sigma)a_\textup{int}(\eta)(1-\td{\td\sigma}) \]
with $\sigma,\td\sigma,\td{\td\sigma}$ as in \eqref{bmu}. This is now a parameter-dependent elliptic element of the cone calculus on $M$ with parameter $\eta\in\R^q$. It is known, see the explanations after this proof, that there is a constant $C>0$ such that the operators \eqref{nnew24} are isomorphisms for all $|\eta|\geq C$. Now, in order to construct $a(\eta)$ such that \eqref{nnew24} are isomorphisms for all $\eta\in\R^q$ we simply perform the construction with $(\eta,\lambda)\in \R^{q+r},r\geq 1$ in place of $\eta$, then obtain a family $a(\eta,\lambda)$ and define $a(\eta):= a(\eta,\lambda^1)$ with a $\lambda^1\in\R^r, |\lambda^1|\geq C$.
\qed \\
Let us now give more information on the above mentioned space
\[ C^\mu(M,\boldsymbol{g};\R^q), \ \boldsymbol{g}=(\gamma,\gamma-\mu), \]
of parameter-dependent cone operators on $M$ of order $\mu\in\R$, with the weight data $\boldsymbol{g}$. The elements $a(\eta)\in C^\mu(M,\boldsymbol{g};\R^q)$ have a principal symbolic hierarchy
\begin{equation} \label{hsmb}
\sigma(a):=(\sigma_\psi(a),\sigma_\wedge(a))
\end{equation}
where $\sigma_\psi(a)$ is the parameter-dependent homogeneous principal symbol of order $\mu$, defined through $a(\eta)\in L^\mu_\clw(M\setminus\{v\};\R^q)$. This determines the reduced symbol
\[ \td\sigma_\psi(a)(t,x,\tau,\xi,\eta):=t^\mu\sigma_\psi(a)(t,x,t^{-1}\tau,\xi,t^{-1}\eta) \]
given close to $v$ in the splitting of variables $(t,x)$ with covariables $(\tau,\xi)$. By construction $\td\sigma_\psi(a)$ is smooth up to $t=0$. The second component $\sigma_\wedge(a)(\eta)$ is defined as 
\begin{multline*} \sigma_\wedge(a)(\eta):=t^{-\mu}\omega_{|\eta|}\sop_M^{\gamma-\frac{n}{2}}(h_0)(\eta) \td\omega_{|\eta|} \\ +t^{-\mu}(1-\omega_{|\eta|})\op_t(p_0)(\eta)(1-\td{\td\omega}_{|\eta|})+ \sigma_\wedge(m+g)(\eta) 
\end{multline*}
where $\sigma_\wedge(m+g)(\eta)$ is just the (twisted) homogeneous principal symbol of $m+g$ as a classical operator-valued symbol.

The element $a(\eta)$ of $C_G(M,\boldsymbol{g};\R^q)$ represent families of continuous operators 
\begin{equation} \label{qetF}
a(\eta):H^{s,\gamma}(M)\to H^{s,\gamma-\mu}(M)
\end{equation}
for all $s\in\R$.
\begin{Def}
An element $a(\eta)\in C^\mu(M,\boldsymbol{g};\R^q)$ is called elliptic, if
\begin{enumerate}
\item[\textup{(i)}] $\sigma_\psi(a)$ never vanishes as a function on $T^*((M\setminus\{v\})\times\R^q)\setminus 0$ and if $\td\sigma_\psi(a)$ does not vanish for all $(t,x,\tau,\xi,\eta)$, $(\tau,\xi,\eta)\not=0$, up to $t=0$;
\item[\textup{(ii)}] $\sigma_\wedge(a)(\eta):\K^{s,\gamma}(X^\wedge)\to\K^{s-\mu,\gamma-\mu}(X^\wedge)$ is a family of isomorphisms for all $\eta\not=0$, and any $s\in\R$. 
\end{enumerate}
\end{Def}
\begin{Thm} \label{htmfred}
If $a(\eta)\in C^\mu(M,\boldsymbol{g};\R^q)$, $\boldsymbol{g}=(\gamma,\gamma-\mu)$ is elliptic, there exists an element $a^{(-1)}(\eta)\in C^{-\mu}(M,\boldsymbol{g}^{-1};\R^q)$ $\boldsymbol{g}^{-1}:=(\gamma-\mu,\gamma)$, such that
\[ 1-a^{(-1)}(\eta)a(\eta)\in C_G(M,\boldsymbol{g}_l;\R^q), \ \ 1-a(\eta)a^{(-1)}(\eta)\in C_G(M,\boldsymbol{g}_r;\R^q), \]
where $\boldsymbol{g}_l:=(\gamma,\gamma)$, $\boldsymbol{g}_r:=(\gamma-\mu,\gamma-\mu)$.  
\end{Thm}
The proof employs known elements of the edge symbolic calculus (cf. \cite{Schu20}); so we do not recall the details here. Let us only note that the inverses of $\sigma_\psi(a),\td\sigma_\psi(a)$ and $\sigma_\wedge(a)$ can be employed to construct an operator family $b(\eta)\in C^{-\mu}(M,\boldsymbol{g}^{-1};\R^q)$ such that 
\[ \sigma_\psi(a)^{(-1)}=\sigma_\psi(b), \ \ \td\sigma_\psi(a)^{(-1)}=\td\sigma_\psi(b), \ \ \sigma_\wedge(a)^{(-1)}=\sigma_\wedge(b). \]
This gives us $1-b(\eta)a(\eta)=:c_0(\eta)\in C^{-1}(M,\boldsymbol{g}_l;\R^q)$, and a formal Neumann series argument allows us to improve $b(\eta)$ to a left parametrix $a^{(-1)}(\eta)$ by setting $a^{(-1)}(\eta):=\left(\sum_{j=0}^\infty c_0^j(\eta)\right) b(\eta)$ (using the existence of the asymptotic sum in $C^0(M,\boldsymbol{g};\R^q)$). In a similar manner we can construct a right parametrix, i.e., $a^{(-1)}(\eta)$ is as desired. 
\begin{Cor} \label{corf1}
If $a(\eta)$ is as in Theorem \textup{\ref{htmfred}}, then \eqref{qetF} is a family of Fredholm operators of index $0$, and there is a constant $C>0$ such that the operators \eqref{qetF} are isomorphisms for all $|\eta|\geq C$, $s\in\R$.
\end{Cor}
\begin{Cor}
If we perform the construction of Theorem \textup{\ref{htmfred}} with the parameter $(\eta,\lambda)\in\R^{q+l}$, $l\geq 1$, rather than $\eta$, Corollary \textup{\ref{corf1}} yields that $a(\eta,\lambda)$ is invertible for all $\eta\in\R^q$, $|\lambda|\geq C$. Then, setting $a(\eta):=a(\eta,\lambda^1)$, $|\lambda^1|\geq C$ fixed, we obtain $a^{-1}(\eta)\in C^{-\mu}(M,\boldsymbol{g}^{-1};\R^q)$.
\end{Cor}
Observe that the operator functions of Theorem \ref{thm132} refer to scales of spaces with two parameters, namely, $s\in\R$, the smoothness, and $\gamma\in\R$, the weight. Compared with Definition \ref{732.symb1901.de} we have here an additional weight. There are two ways to make the different view points compatible. One is to apply weight reducing isomorphisms 
\begin{equation} \label{h-gamm}
\textup{h}^{-\gamma}:H^{s,\gamma}(M)\to H^{s,\gamma-\mu}(M)
\end{equation}
in \eqref{equiso}. Then, passing from 
\begin{equation} \label{new18}
a(\eta):H^{s,\gamma}(M)\to H^{s-\mu,\gamma-\mu}(M)
\end{equation}
to 
\begin{equation} \label{new1,18}
b^\mu(\eta):=\textup{h}^{-\gamma+\mu}a(\eta)\textup{h}^\gamma:H^{s,0}(M)\to H^{s-\mu,0}(M)
\end{equation}
we obtain operator functions between spaces only referring to $s$ but with properties as required in Definition \ref{732.symb1901.de} (which remains to be verified).
\begin{Rem}
The spaces $E^s:=H^{s,0}(M),\, s\in\R$, form a scale with the properties at the beginning of Section \textup{\ref{OrdRed}}.
\end{Rem}
Another way is to modify the abstract framework by admitting scales $E^{s,\gamma}$ rather than $E^s$, where in general $\gamma$ may be in $\R^k$ (which is motivated by the higher corner calculus). We do not study the second possibility here but we only note that the variant with $E^{s,\gamma}$-spaces is very similar to the one without $\gamma$. \\
Let us now look at operator functions of the form \eqref{new1,18}.
\begin{Thm}
\label{theormb}
The operators \eqref{new1,18} constitute an order reducing family in the spaces $E^s:=H^{s,0}(M)$, where the properties \textup{(i)-(iii)} of Definition \textup{\ref{ordred}} are satisfied. 
\end{Thm}
\proof 
In this proof we concentrate on the properties of our operators for every fixed $s,\mu,\nu$ with $\nu\geq\mu$. The uniformity of the involved constants can easily be deduced; however, the simple (but lengthy) considerations will be left out.\\
(i) We have to show that 
\[ D_\eta^\beta b^\mu(\eta)=D_\eta^\beta\{\textup{h}^{-\gamma+\mu} a(\eta)\textup{h}^\gamma\} \in C^\infty(\R^q,\Li(E^s,E^{s-\mu+|\beta|})) \]
for all $s\in\R$, $\beta\in\N^q$. According to \eqref{bmu} the operator function is a sum of two contributions. The second summand
\[ (1-\sigma)\textup{h}^{-\gamma+\mu}a_\textup{int}(\eta)\textup{h}^\gamma(1-\td{\td\sigma}) \]
is a parameter-dependent family in $L^\mu_\clw(2\M;\R^q)$ and obviously has the desired property. The first summand is of the form 
\[\sigma\textup{h}^{-\gamma+\mu}\{a_\textup{edge}(\eta)+m(\eta)+g(\eta)\}\textup{h}^\gamma\td\sigma.\]
From the proof of  Theorem \ref{thm132} we have
\[ D_\eta^\beta \sigma  a_\textup{edge}(\eta)\td\sigma\in S^{\mu-|\beta|}(\R^q;\K^{s,\gamma;g} (X^\wedge),\K^{s-\mu+|\beta|,\gamma-\mu;g}(X^\wedge))\]
for every $\beta\in\N^q$. In particular, these operator functions are smooth in $\eta$ and the derivates improve the smoothness in the image by $|\beta|$. This gives us the desired property of $\sigma \textup{h}^{-\gamma+\mu}a_\textup{edge}(\eta)\textup{h}^\gamma\td\sigma$. The $C^\infty$ dependence of $m(\eta)+g(\eta)$ in $\eta$ is clear (those are operator-valued symbols), and they map to $\K^{\infty,\gamma-\mu;g}(X^\wedge)$ anyway. Therefore, the desired property of $\sigma \textup{h} ^{-\gamma+\mu}\{m(\eta)+g(\eta)\}\textup{h}^\gamma\td\sigma$ is satisfied as well.\\
(ii) This property essentially corresponds to the fact that the product in consideration close to the conical point is a symbol in $\eta$ of order zero and that the group action in $\K^{0,0}(X^\wedge)$-spaces is unitary. Outside the conical point the boundedness is as in Example \ref{ExamX}.\\
(iii) The proof of this property close to the conical point is of a similar structure as Proposition \ref{nn9}, since our operators are based on operator-valued symbols referring to spaces with group action. The contribution outside the conical point is as in Example \ref{ExamX}.
\qed
\begin{Rem}
For $E^s:=\Hi^{s,0}(M)$, $s\in\R$, $\E=(E^s)_{s\in\R}$, the operator functions $b^\mu(\eta)$ of the form \eqref{new1,18} belong to $S^\mu (\R^q;\E,\E)$ \lr see the notation after Definition \textup{\ref{732.symb1901.de}}\rr.
\end{Rem}

\section{Operators referring to a corner point}
\label{Oprefex}
\subsection{Weighted spaces}

Let ${\E} = (E^s)_{s \in \R} \in \mathfrak{E}$ be a scale
and $(b^\mu(\varrho))_{\mu \in \R}$, $\varrho \in \R$, 
be an order reducing
family (see Definition \ref{ordred} with 
$q = 1$). We 
define a new scale of spaces adapted to the Mellin
transform and the approach of the cone calculus.
In the following definition the Mellin transform refers
to the variable $r \in \R_+$, i.e., $M = M_{r \to w}$.

\begin{Def}
\label{733.1901.de}
For every $s, \gamma \in \R$ we define the space $\Hi^{s, \gamma}(\R_+,\E)$ to be the completion of
$C_0^\infty(\R_+, E^\infty)$ with respect to be norm
\begin{equation}
\label{733.208neu2201.eq}  
\| u \|_{ {\Hi}^{s, \gamma}(\R_+, {\E})} =
      \Bigl\{ \frac{1}{2 \pi i}
      \int_{\Gamma_{\frac{d+1}{2} - \gamma}} \| 
      b^s(\imb w)
      (Mu)(w) \|_{E^0}^2 dw \Bigr\}^{\frac{1}{2}}   
\end{equation}
for a $d = d_{\cal E} \in \N$.
The Mellin transform $M$ in \eqref{733.208neu2201.eq} is
interpreted as the weighted Mellin transform $M_{\gamma -
\frac{d}{2}}$.
\end{Def}      

The role of $d_{\E}$ is an extra information, given
together with the scale ${\E}$. In the example ${\E} =
(H^s(X))_{s \in \R}$ for a closed compact $C^\infty$
manifold $X$ we set $d_{\E} := \dim X$.

Observe that when we replace the order reducing family in
\eqref{733.208neu2201.eq} by an equivalent one the
resulting norm is equivalent to
\eqref{733.208neu2201.eq}.

By virtue of the identity
\[  r^\beta {\Hi}^{s, \gamma}(\R_+, {\E}) = {\Hi}^{s, \gamma + \beta}(\R_+, {\E})    \]
for every $s, \gamma, \beta \in \R$, it is often enough to refer the considerations to one
particular weight, or to set
\begin{equation}
\label{733.(d)2201.eq}
d_{\E} = 0.    
\end{equation}
For simplicity we now assume 
\eqref{733.(d)2201.eq}.

Let us consider Definition \ref{732.symb1901.de} for the case
$U = \R$, $q = 1$, and denote the covariable now by
$\varrho \in \R$. Set
\[  S^\mu(\Rb_+ \times \Rb_+ \times \R; {\E},
     \widetilde{\E}) := S^\mu(\R \times \R \times \R;
     {\E}, \widetilde{\E}) |_{\Rb_+ \times \Rb_+ \times \R} \]
and
\begin{align*}  
S^\mu (\Rb_+ \times \Rb_+ \times \Gamma_\delta;
     \E, \widetilde{\E}) := \ & \{ a(r,r',w) \in
     C^\infty(\Rb_+ \times \Rb_+ \times
     \Gamma_\delta, 
      {\Li}^\mu({\E}, \widetilde{\E})) \\
      & : a(r, r', \delta + i \varrho) \in
       S^\mu(\Rb_+ \times \Rb_+ \times \R_\rho; {\E}, \widetilde{\E}) \}   
\end{align*}       
for any $\delta \in \R$. The subspaces of
$r'$-independent ($(r,r')$-independent) symbols are
denoted by $S^\mu(\Rb_+ \times \R; {\E},
\widetilde{\E})$ \ $(S^\mu(\R; {\E},
\widetilde{\E}'))$ and
$S^\mu({\Rb}_+ \times \Gamma_\delta; {\E},
\widetilde{\E})$ \ $(S^\mu(\Gamma_\delta; {\E},
\widetilde{\E}))$, respectively.

Given an element $f(r,r',w) \in S^\mu({\Rb}_+ \times
{\Rb}_+ \times \Gamma_{\frac{1}{2} - \gamma}; {\E},
\widetilde{\E})$ we set             
\begin{equation}
\label{733.Me2201.eq}
\sop_M^\gamma(f)u(r) = \frac{1}{2 \pi} \int
   \int_{0}^{\infty} (\frac{r}{r'})^{-(\frac{1}{2} - 
   \gamma + i \varrho)} f(r,r', \frac{1}{2} - \gamma + i
   \varrho) u (r') \frac{dr'}{r'} d \varrho.
\end{equation}
Let, for instance, $f$ be independent of $r'$. Then
\eqref{733.Me2201.eq} induces a continuous operator
\begin{equation}
\label{733.Mo2201.eq}
\sop_M^\gamma(f) : C_0^\infty(\R_+, E^\infty) \to
     C^\infty(\R_+, \widetilde{E}^\infty).    
\end{equation}  

In fact, we have $\sop_M^\gamma(f) = M_{\gamma,w \to
r}^{-1} f(r,w) M_{\gamma,r' \to w}$. The weighted Mellin 
transform $M_\gamma$ induces a continuous operator
\[  M_\gamma : C_0^\infty(\R_+, E^s) \to {\Sch}(\Gamma_{\frac{1}{2} - \gamma}, E^s)    \]   
for every $s \in \R$. The subsequent multiplication of
$M_\gamma u(w)$ by $f(r,w)$ gives rise to an element in
$C^\infty(\R_+, {\Sch}(\Gamma_{\frac{1}{2} - \gamma},
\widetilde{E}^{s- \mu}))$, and then it follows
easily that $\sop_M^\gamma(f) u \in C^\infty(\R_+,
\widetilde{E}^{s- \mu})$. We now formulate a
continuity result, first for the case of symbols with
constant coefficients.

\begin{Thm} \label{Thm31227}
\label{733.cont2201.th}
For every $f(w) \in S^\mu(\Gamma_{\frac{1}{2} - \gamma};
{\E}, \widetilde{\E})$ the operator
\eqref{733.Mo2201.eq} extends to a continuous operator
\begin{equation}
\label{733.cof2201.eq}
\sop_M^\gamma(f) : {\Hi}^{s, \gamma}(\R_+, {\E})
     \to {\Hi}^{s- \mu,\gamma}(\R_+, \widetilde{\E})
\end{equation}
for every $s \in \R$. Moreover, $f \to \sop_M^\gamma(f)$
induces a continuous operator
\begin{equation}
\label{733.nof2201.eq}
S^\mu(\Gamma_{\frac{1}{2} - \gamma}; {\E},
    \widetilde{\E}) \to {\Li}({\Hi}^{s,
    \gamma}(\R_+, {\E}), {\Hi}^{s- \mu,
    \gamma}(\R_+, \widetilde{\E}))
\end{equation}
for every $s \in \R$.
\end{Thm}
\proof
We have
\begin{align*}
\| \sop_M^\gamma(f) u & \|_{ {\Hi}^{s-\mu,\gamma}(\R_+,
        \widetilde{\E})}^2 \\
  = & \int_\R
      \| \tilde{b}^{s- \mu}(\varrho) M_\gamma
       (M_\gamma^{-1} 
        f(\frac{1}{2} - \gamma + i \varrho)) (M_\gamma u)
	 (\frac{1}{2} - \gamma + i \varrho)
	 \|_{\widetilde{E}^0}^2 d \varrho \\	
  = & \int_\R \|
        \tilde{b}^{s- \mu}(\varrho) f(\frac{1}{2} - 
	\gamma + i \varrho)b^{- s}(\varrho) b^s(\varrho)
	(M_\gamma u) (\frac{1}{2} - \gamma + i \varrho)
	\|_{\widetilde{E}^0}^2 d \varrho \\
	\leq & \,\, c^2 \| u
	\|^2_{ { \Hi}^{s, \gamma}(\R_+, {\E})} 
\end{align*}
with	    
\[  c = \sup_{\varrho \in \R} \| \tilde{b}^{s-
      \mu}(\varrho) f(\frac{1}{2}- \gamma + i
      \varrho)b^{-s}(\varrho) \|_{ {\Li}(E^0,\widetilde{E}^0)}     \]         
which is finite for every $s \in \R$ (cf. the estimates
\eqref{732.est1901.eq}). Thus we have proved the
continuity both of \eqref{733.cof2201.eq} and
\eqref{733.nof2201.eq}.
\qed 

In order to generalise Theorem \ref{Thm31227} to symbols with variable coefficients we impose conditions of reasonable generality that allow us to reduce the arguments to a vector-valued analogue of Kumano-go's technique.\\
Given a Fr\'echet space $V$ with a countable system of semi-norms $(\pi_\iota)_{\iota\in\N}$ that defines its topology, we denote by
\[ C_B^\infty (\R_+\times \R_+ ,V) \]
the set of all $u(r,r')\in C^\infty(\R_+\times\R_+,V)$ such that 
\[ \sup_{r,r'\in\R_+} \pi_\iota\Big((r\partial_r)^k(r'\partial_{r'})^{k'}u(r,r')\Big) < \infty \]
for all $k,k'\in\N$. In a similar manner by $C_B^\infty (\R_+ ,V)$ we denote the set of such functions that are independent of $r'$.\\
Moreover, we set 
\[ S_B^\mu(\R_+\times\R_+\times\Gamma_{\oneh-\gamma};\E,\widetilde\E):= C_B^\infty (\R_+\times \R_+, S^\mu(\Gamma_{\oneh-\gamma};\E,\widetilde\E)) \]
and, similarly, $S_B^\mu(\R_+\times\Gamma_{\oneh-\gamma};\E,\widetilde\E):= C_B^\infty(\R_+,S^\mu(\Gamma_{\oneh-\gamma};\E, \widetilde\E))$.
\begin{Thm} \label{ThmconT}
For every $f(r,w)\in S_B^\mu(\R_+\times\Gamma_{\oneh-\gamma};\E,\widetilde\E)$ the operator $\sop_M^\gamma(f)$ induces a continuous mapping
\[ \sop_M^\gamma(f):\Hi^{s,\gamma}(\R_+,\E)\to\Hi^{s-\mu,\gamma}(\R_+,\widetilde\E) , \]
and $f\to \sop_M^\gamma(f)$ a continuous operator 
\[ S_B^\mu(\R_+\times\Gamma_{\oneh-\gamma};\E,\widetilde\E)\to \Li(\Hi^{s,\gamma}(\R_+,\E), \Hi^{s-\mu,\gamma}(\R_+,\widetilde\E)) \]
for every $s\in\R$.
\end{Thm}

Parallel to the spaces of Definition \ref{733.1901.de} it
also makes sense to consider their ``cylindrical''
analogue, defined as follows.

\begin{Def}
\label{733.HS2201.de}
Let $(b^s(\eta))_{s \in \R}$, be an order reducing family
as in Definition \textup{\ref{ordred}}. For every $s \in
\R$ we define the space $H^s(\R^q,\E)$ to be the
completion of $C_0^\infty(\R^q, E^\infty)$ with respect
to the norm
\[  \| u \|_{H^s(\R^q,\E )} := \Bigl\{\int_{\R^q} \|
       b^s(\eta) (\Fu u)(\eta) \|_{E^0}^2 d \eta 
       \Bigr\}^{\frac{1}{2}}.   \]
\end{Def}

Clearly, similarly as above, with a symbol $a(y,y',\eta)
\in S^\mu(\R^q \times \R^q \times \R^q; {\cal E},
\widetilde{\E})$ (when we impose a suitable 
control with respect to
the dependence on $y'$ for large $|y'|$) we can associate
a pseudo-differential operator
\[  \op_y(a) u(y) = \iint e^{i(y-y')\eta}a(y,y',\eta)
      u(y') dy' \dbar \eta.    \]
In particular, if $a = a(\eta)$ has constant
coefficients, then we obtain a continuous operator
\[  \op_y(a) : H^s(\R^q, {\E}) \to H^{s- \mu}(\R^q,
      \widetilde{\E})     \]
for every $s \in \R$. In the case of variable
coefficients we need some precautions on the nature of
symbols. This will be postponed for the moment.   

We are mainly interested in the case $q = 1$. Consider
the transformation
\[  (S_\gamma u)(y) := e^{- (\frac{1}{2} - \gamma)y}
       u(e^{- y})      \]                  
from functions in $r \in \R_+$ to functions in $y \in \R$. We then have the identity
\[  (M_\gamma u)(\frac{1}{2} - \gamma + i \varrho) = (\Fu S_\gamma u)(\varrho)    \]
with $\Fu$ being the one-dimensional Fourier transform.
This gives us
\[  \Bigl\{ \frac{1}{2 \pi} \int_{\R} \| b^s(\eta) (\Fu S_\gamma u)(\eta) \|_{E^0}^2 d \eta 
	\Bigr\}^{\frac{1}{2}} = \| S_\gamma u \|_{H^s(\R,
	{\E})} = \| u \|_{ {\Hi}^{s,\gamma}(\R_+,\E)},            \]
i.e., $S_\gamma$ induces an isomorphism
\[  S_\gamma : {\Hi}^{s, \gamma}(\R_+, {\E}) \to H^s(\R, {\E}).   \]

\begin{Rem}
\label{733.re2201.re}
By reformulating the expression \eqref{733.Me2201.eq} we
obtain
\[ \sop_M^\gamma(f) u(r) = \frac{1}{2 \pi} \iint
      e^{(\frac{1}{2}-\gamma+i\varrho)(\log r' - \log r)}
      f(r,r', \frac{1}{2} - \gamma + i \varrho)
      u(r') \frac{dr'}{r'} d \varrho.    \]
Substituting $r = e^{- y}$, $r' = e^{-y'}$ gives us
\begin{align*}  
\sop_M^\gamma(f) u(r) = & \frac{1}{2 \pi} \iint
       e^{i(y-y')\varrho} e^{( \frac{1}{2} - \gamma)
          (y-y')} f(e^{- y}, \\ 
       & e^{- y'}, \frac{1}{2} - \gamma + i
       \varrho) u(e^{- y'}) dy' d \varrho   
       = \op_y(g_\gamma)v(y)     
\end{align*}                         	
with
$v(y) := u(e^{- y})$
and $g_\gamma(y,y', \varrho) := e^{(\frac{1}{2} -
      \gamma)(y - y')} f(e^{- y}, e^{- y'}, \frac{1}{2} -
      \gamma + i \varrho)$.    

In other words, if $\chi: \R_+ \to \R$ is defined by
$\chi(r) = - \log r =: y$, we have $(\chi^* v)(r) = v(-
\log r)$ or $((\chi^{-1})^* u)(g) = u(e^{- y})$ and
\[  \sop_M^\gamma(f) = \chi^*
        \op_y(g_\gamma)(\chi^{-1})^*.     \] 
Thus $\op_y(g_\gamma)$ is the operator push forward of
$\sop_M^\gamma(f)$ under $\chi$.
\end{Rem}	         

\subsection{Mellin quantisation and kernel cut-off}

The axiomatic cone calculus that we develop here is a substructure of the general calculus of operators with symbols in $a(r,\rho) \in S^\mu(\R_+\times\R;\E,\widetilde\E)$ of the form $a(r,\rho)=\td a(r,r\rho), \td a(r,\td\rho)\in S^\mu(\Rb_+\times\R_{\td\rho}; \E,\widetilde\E)$ (up to a weight factor and modulo smoothing operators) with a special control near $r=0$ via Mellin quantisation. By $L^{-\infty}(\R_+;\E,\widetilde\E;\R^q)$ we denote the space of all Schwartz functions in $\eta\in\R^q$ with values in operators
\[ C_0^\infty(\R_+,E^{-\infty})\to C^\infty(\R_+,\widetilde E^\infty). \]
We then define
\begin{multline}
L^\mu(\R_+;\E,\widetilde\E;\R^q)=\{\op_r(a)(\eta)+C(\eta): \\
\nonumber a(r,\rho,\eta)\in S^\mu(\R_+\times\R^{1+q}_{\rho,\eta};\E, \widetilde\E), C(\eta)\in L^{-\infty}(\R_+;\E,\widetilde\E;\R^q)\}. 
\end{multline}
Our next objective is to formulate a Mellin quantisation result of symbols 
\begin{equation}
\label{quantsym}
a(r,\rho,\eta)= \td a(r,r\rho,r\eta),\ \ \td a(r,\td\rho,\td\eta)\in S^\mu (\Rb_+\times\R^{1+q}_{\td\rho,\td\eta}; \E,\widetilde\E)
\end{equation}
(see Remark \ref{Frechnat}).
\begin{Def} \label{defmo}
By $M_\mathcal{O}^\mu(\E,\widetilde\E;\R^q_{\td\eta})$ we denote the set of all $h(z,\td\eta)\in \mathcal{A}(\C , S^\mu(\R^q_{\td\eta};\E,\widetilde\E))$ such that 
\[ h(\beta+i\rho,\td\eta)\in S^\mu (\R_\rho\times\R^q_{\td\eta} ;\E,\widetilde\E) \]
for every $\beta\in\R$, uniformly in compact $\beta$-intervals. For $q=0$ we simply write $M_\mathcal{O}^\mu(\E,\widetilde\E)$.
\end{Def}
\begin{Thm}\label{melquan2}
For every symbol $a(r,\rho,\eta)$ of the form \eqref{quantsym} there exists an $\td h(r,z,\td\eta)\in C^\infty (\Rb_+,M_\mathcal{O}^\mu(\E,\widetilde\E;\R^q))$ such that for $h(r,z,\eta):=\td h (r,z,r\eta)$ and every $\delta\in\R$ we have 
\[ \sop_M^\delta (h)(\eta)=\op_r(a)(\eta)  \]
modulo operators in $L^{-\infty}(\R_+;\E,\widetilde\E;\R^q)$.
\end{Thm}
This result in the context of operator-valued symbols based on order reductions is mentioned here for completeness. It is contained in a joint work (in preparation) of the second author with C.-I. Martin (Potsdam) and N. Rablou (G\"ottingen). It extends a corresponding result of the edge symbolic calculus, see \cite[Theorem 3.2]{Gil2}. More information in that case is given in \cite[Chapter 4]{Krai9}. Here we adapt some part of this approch to realise the kernel cut-off principle that allows us to recognise how many parameter-dependent meromorphic Mellin symbols exist.
\begin{Def}
Let $S^\mu(\boldsymbol{C}\times\R^q;\E,\widetilde\E)$ denote the space of all $h(\zeta,\eta)\in\mathcal{A}(\C,S^\mu(\R^q_\eta;\E,\widetilde\E))$ such that 
\[ h(\rho+i\delta,\eta)\in S^\mu(\R^{1+q}_{\rho,\eta};\E,\widetilde\E) \]
for every $\delta\in\R$, uniformly in compact $\delta$-intervals.
\end{Def}
Clearly the space $S^\mu(\boldsymbol{C}\times\R^q;\E,\widetilde\E)$ is a generalisation of $M^\mu_\mathcal{O}(\E,\widetilde\E)$, however, with an interchanged role of real and imaginary part of the complex covariable. To produce elements of $S^\mu(\boldsymbol{C}\times\R^q;\E,\widetilde\E)$ we consider a so-called kernel cut-off operator
\[ V:C_0^\infty(\R) \times S^\mu(\R^{1+q};\E,\widetilde\E)\to S^\mu(\boldsymbol{C}\times\R^q;\E,\widetilde\E)\]
transforming an arbitrary element $a(\rho,\eta)\in S^\mu(\R^{1+q};\E,\widetilde\E)$ into $\big( V(\varphi)a\big)(\zeta,\eta)\in S^\mu(\boldsymbol{C}\times\R^q;\E,\widetilde\E)$ for any $\varphi \in C_0^\infty(\R)$. It will be useful to admit $\varphi$ to belong to the space $C^\infty_\textup{b}(\R):=
\{\varphi\in C^\infty(\R_\theta):\sup_{\theta\in\R}|D^k_\theta \varphi(\theta)|<\infty\ \textup{for every} \ k\in\N\}$.
We set
\begin{equation} \label{oscorn}
\big(V(\varphi)a\big)(\rho,\eta):=\iint e^{-i\theta\td\rho}\varphi(\theta)a(\rho-\td\rho,\eta)d\theta 
\dbar\td\rho ,
\end{equation}
interpreted as an oscillatory integral (see also \cite{Kuma1}). We now prove the following result:
\begin{Thm}
\label{734.cutR0703.th}
The kernel cut-off operator
\[  V : (\varphi, a) \to V(\varphi)a   \]
defines a bilinear and continuous mapping
\begin{equation}
\label{734.cont0703.eq}  
V : C^\infty_\textup{b}(\R) \times S^\mu(\R^{1+q}; {\cal E},
       \widetilde{\cal E}) \to S^\mu(\R^{1+q}; {\cal E},
       \widetilde{\cal E}),
\end{equation}
and $\big(V(\varphi) a\big)(\rho, \eta)$ admits an asymptotic
expansion
\begin{equation}
\label{734.asV0703.eq}
\big(V(\varphi)a\big)(\rho, \eta) \sim \sum_{k=0}^{\infty}
\frac{(-1)^k}{k!} D_\theta^k \varphi(0)
\partial_\rho^k a(\rho, \eta).
\end{equation}
\end{Thm}
\proof
First note that the mapping
\begin{align*}
(\varphi, a) & \to \varphi(\theta) a(\rho - \tilde{\rho}, \eta),       \\
& C^\infty_\textup{b}(\R) \times S^\mu(\R_{\rho, \eta}^{1+q};
{\cal E}, \widetilde{\cal E}) \to C^\infty(\R_{\rho,
    \eta}^q, S^\mu_\textup{b}(\R_\theta \times \R_{\tilde{\rho}};
    {\cal E}, \widetilde{\cal E}))
\end{align*}
for $S^\mu_\textup{b}(\R_\theta \times \R_{\tilde{\rho}}; {\cal E}, 
\widetilde{\cal E}) := C^\infty_\textup{b}(\R_\theta,
S^\mu(\R_{\tilde{\rho}}; {\cal E}, \widetilde{\cal
E}))$ is bilinear and continuous. For the proof of the
continuity of \eqref{734.cont0703.eq} it suffices to
verify that $\big(V(\varphi)a\big)(\rho, \eta) \in
S^\mu(\R^{1+q}; {\cal E}, \widetilde{\cal E})$ and then
to apply the closed graph theorem. By virtue of
\[D_{\rho, \eta}^\beta \big(V(\varphi) a\big)(\rho, 
\eta) = \big(V(\varphi)(D_{\rho,
\eta}^\beta a)\big)(\rho, \eta) \]
for every $\beta \in
\N^{1+q}$ we only have to check that for every $s \in \R$
\begin{equation}
\label{734.(P)0703.eq}
\| \tilde{b}^{s-\mu}(\rho, \eta) \big(V(\varphi) a\big)(\rho,
   \eta) b^{-s}(\rho, \eta) \|_{ {\cal L}(E^0,
   \widetilde{E}^0)} \leq c
\end{equation}
for all $(\rho, \eta) \in \R^{1+q}$, with a constant
$c = c(s) > 0$. We regularise the oscillatory integral
\eqref{oscorn}
\[  \big(V(\varphi) a\big)(\rho, \eta) = \iint e^{-i \theta
      \tilde{\rho}}\langle \theta \rangle^{-2} \{ (1-
      \partial_\theta^2)^N \varphi(\theta) \}
      a_N(\rho, \tilde{\rho},\eta) d \theta \dbar
      \tilde{\rho}    \]
for
\begin{equation}
\label{734.est0703.eq}
a_N(\rho, \tilde{\rho}, \eta) := (1-
    \partial_{\tilde{\rho}}^2) \{ \la
    \tilde{\rho} \ra^{-2 N} a(\rho -
    \tilde{\rho}, \eta) \}.
\end{equation}
 The function \eqref{734.est0703.eq} is a linear
 combination of terms
 \[  (\partial_{\tilde{\rho}}^j \langle
       \tilde{\rho} \rangle^{- 2 N}) 
       (\partial_\rho^k a)
       (\rho - \tilde{\rho}, \eta) \ \textup{for $0
       \leq j,k \leq 2$}.    \]

We have                   
\begin{align}
\label{734.estb0803.eq}
\begin{aligned}
\| \tilde{b}^{s-\mu}&(\rho, \eta) \Bigl\{ \iint e^{-i
      \theta \tilde{\rho}} 
       \langle \theta \rangle^{-2}
       (1-\partial_\theta^2)^N \varphi(\theta)
      (\partial_{\tilde{\rho}}^j \langle
      \tilde{\rho} \rangle^{-2N}) \\
   & \ (\partial_\rho^k
      a) (\rho-\tilde{\rho}, \eta) d \theta \dbar
      \tilde{\rho} \Bigr\} b^{-s}(\rho, \eta) 
      \|_{ {\cal L}(E^0, \widetilde{E}^0)}	\\		
 = & \ \| \iint \tilde{b}^{s-\mu}(\rho, \eta)
      \tilde{b}^{-s+\mu}(\rho- \tilde{\rho}, \eta)
      \tilde{b}^{s- \mu}(\rho- \tilde{\rho}, \eta)      		
    \bigl\{ e^{-i\theta \tilde{\rho}} \langle
      \theta \rangle^{-2}  
            (1- \partial_\theta^2)^N  \varphi(\theta)  \\                      
&   (\partial_{\tilde{\rho}}^j
       \langle \tilde{\rho} \rangle^{-2N}) 
    (\partial_\rho^k a) (\rho- \tilde{\rho},
        \eta) \} b^{-s}(\rho - \tilde{\rho}, \eta)
        b^s(\rho- \tilde{\rho}, \eta)
      b^{-s}(\rho, \eta) d \theta \dbar \tilde{\rho}
        \|_{ {\cal L}(E^0, E^0)} \\
\leq & \ c \iint \| \tilde{b}^{s-\mu} (\rho, \eta)
      \tilde{b}^{-s+\mu}(\rho-\tilde{\rho}, 
      \eta) \|_{  {\cal L}(\widetilde{E}^0, 
      \widetilde{E}^0)}
        \| \tilde{b}^{s-\mu}(\rho- \tilde{\rho},
	\eta) (\partial_{\tilde{\rho}}^j \langle
	\tilde{\rho} \rangle^{-2N}) \\
  &	(\partial_\rho^k
	a) (\rho- \tilde{\rho}, \eta)
	b^{-s}(\rho- \tilde{\rho}, \eta) \|_{ {\cal
	L}(E^0, \widetilde E^0)}    
     \ \| b^s(\rho - \tilde{\rho}, \eta) b^{-s}
          (\rho, \eta) \|_{ {\cal L}(E^0, E^0)} \dbar
	  \tilde{\rho}.
\end{aligned}
\end{align}	  	

For the norms under the integral we apply the Taylor
expansion

\begin{align*}  
b^s(\rho-\tilde{\rho}, \eta) = & \sum_{m=0}^{M}
    \frac{1}{m!} (\partial_\rho^m b^s)(\rho, \eta)
    (- \tilde{\rho})^m   \\
& + \frac{\tilde{\rho}}{M!}^{M+1} \int_{0}^{1} (1-t)^M
      (\partial_\rho^{M+1} b^s)(\rho - t
      \tilde{\rho},\eta) dt.
\end{align*}
This yields          
\begin{align*}
\| b^s & (\rho-\tilde{\rho}, \eta) b^{-s}(\rho,
          \eta) \|_{ {\cal L}(E^0,E^0)}  \\
\leq & \sum_{m=0}^{M} \frac{1}{m!} \langle
       \tilde{\rho} \rangle^m \| (\partial_\rho^m
       b^s) (\rho, \eta) b^{-s} (\rho, \eta) \|_{
       {\cal L}(E^0, E^0)} \\
 & + \frac{\langle \tilde{\rho} \rangle^{M+1}}{M!}
       \int_{0}^{1} (1- t)^M \| (\partial_\rho^{M+1}
       b^s)(\rho - t \tilde{\rho}, \eta) b^{-s} 
       (\rho, \eta) \|_{ {\cal L}(E^0, E^0)} dt.
\end{align*}               	  

By virtue of \eqref{ordsmb}, Proposition
\ref{732.2201.pr} and Proposition \ref{732.a2201.pr} we
obtain $\|(\partial_\rho^m b^s)(\rho,
\eta)b^{-s}(\rho, \eta) \|_{ {\cal L}(E^0, E^0)} \leq
c \langle \rho, \eta \rangle^{-m}$.
Moreover, using Definition \ref{ordred} (iii), it follows that
\begin{align*}
\|(\partial_\rho^{M+1} b^s)&(\rho - t
     \tilde{\rho}, \eta) b^{-s} (\rho - t
     \tilde{\rho}, \eta) b^s(\rho- t
     \tilde{\rho}, \eta) b^{-s}(\rho, \eta) \|_{
     {\cal L}(E^0, E^0)}  \\
\leq & \ c \| (\partial_\rho^{M+1} b^s)(\rho - t
         \tilde{\rho}, \eta) b^{-s} (\rho - t
         \tilde{\rho}, \eta) \| _{ {\cal L}(E^0, E^0)}
	 \\
     & \| b^s(\rho - t \tilde{\rho}, \eta) 
          \|_{ {\cal
            L}(E^s, E^0)} \| b^{-s} (\rho,
	    \eta) \|_{ {\cal L}(E^0, E^s)} \\
 \leq & \ \langle \rho - t \tilde{\rho}, \eta
          \rangle^{B_1(s)} \langle \rho, \eta
	  \rangle^{B_2(s)}
\end{align*}
with certain $B_i(s)$, $i = 1,2$. We thus obtain	     	    	  
\begin{align*}
\| b^s& (\rho - \tilde{\rho}, \eta) b^{-s} (\rho,
       \eta) \|_{ {\cal L}(E^0, E^0)}   \\
 \leq & \ c \langle \tilde{\rho} \rangle^{M+1}
         (\sup_{|t|\leq 1} \langle \rho - t
	 \tilde{\rho}, \eta 
	 \rangle^{-(M+1) + B_1(s)}) 
	 \langle \rho, \eta \rangle^{B_2(s)}.
\end{align*}
By Peetre's inequality for $L \geq 0$ we have
\[  \sup_{|t|\leq 1} \langle \rho - t \tilde{\rho},
       \eta \rangle^{- L} \leq c \langle \tilde{\rho} 
       \rangle^L \langle \rho, \eta \rangle^{- L}. \]	 
Thus choosing $M$ so large that
\[  - (M+1) + B_1(s) \leq 0, \
       - (M+1) + B_1(s) + B_2(s) \leq 0,     \]
it follows that
\begin{align}
\label{734.256neu1303.eq}
\| b^s&(\rho - \tilde{\rho}, \eta) b^{-s} (\rho,
     \eta) \|_{ {\cal L}(E^0, E^0)} \\
  \leq  & \ c \langle \tilde{\rho} \rangle^{M+1}
           \langle \tilde{\rho} \rangle^{M+1- B_1(s)}
	   \langle \rho, \eta \rangle^{-(M+1) + B_1(s)
	   + B_2(s)}
	 \leq c \langle \tilde{\rho} \rangle^{A(s)}           
                       \nonumber
\end{align}
for $A(s) := 2(M+1) - B_2(s)$.

In a similar manner we can show that
\begin{equation}
\label{734.(A2)0803.eq}
\| \tilde{b}^{s- \mu}(\rho, \eta)
     \tilde{b}^{-s+\mu}(\rho - \tilde{\rho}, \eta)
     \|_{ {\cal L}(\widetilde{E}^0, \widetilde{E}^0)}
     \leq c \langle \tilde{\rho}
     \rangle^{\widetilde{A}(s)}    
\end{equation}
for some $\widetilde{A}(s) \in \R$. Applying 
\eqref{734.256neu1303.eq} and \eqref{734.(A2)0803.eq} 
in the estimate \eqref{734.estb0803.eq} it follows that
\begin{equation}
\label{734.(E)0803.eq}
\| \tilde{b}^{s-\mu} (\rho, \eta) \big(V(\varphi)
      a\big)(\rho, \eta) b^{-s} (\rho, \eta) \|_{ {\cal
      L}(E^0, \widetilde{E}^0)}     
      \leq c \sum_{0 \leq j \leq 2} \int
      |\partial_{\tilde{\rho}}^j \langle
      \tilde{\rho} \rangle^{-2N} | \langle
      \tilde{\rho} \rangle^{A(s) + \widetilde{A}(s)} 
      \dbar \tilde{\rho}.
\end{equation}
Since  $N \in \N$ can be chosen as large as we want, it
follows that the right hand side of
\eqref{734.(E)0803.eq} is finite for an appropriate $N$.
This completes the proof of \eqref{734.(P)0703.eq}.
The relation \eqref{734.asV0703.eq} immediately follows by applying the Taylor expansion of $\vp$ at $0$.  
\qed \\
\begin{Thm}
\label{734.hol0803.th}
The kernel cut-off operator $V : (\varphi, a) \to
V(\varphi) a$ defines a bilinear and continuous mapping
\begin{equation}
\label{734.conthol803.eq}
V : C_0^\infty(\R) \times S^\mu(\R^{1+q}; {\cal E},
      \widetilde{\cal E}) \to S^\mu(\boldsymbol{C} \times \R^q;
      {\cal E}, \widetilde{\cal E}).
\end{equation}
\end{Thm}    
\proof
Writing
\[  \big(V(\varphi)a\big)(\rho, \eta) = \int e^{-i \theta
      \rho} \varphi(\theta) \bigl\{ \int e^{i \theta \rho'}  
      a(\rho', \eta) \dbar \rho' \bigr\} 
      d \theta     \]
we see that $\big(V(\varphi)a\big)(\rho, \eta)$ is the Fourier
transform of a distribution $\varphi(\theta) \int
e^{i \theta \rho'} a(\rho', \eta) \dbar \rho'
\in {\cal S}'(\R_\theta, {\cal L}^\mu({\cal E},
\widetilde{\cal E}))$ with compact support. This extends
to  a holomorphic ${\cal L}^\mu({\cal E}, \widetilde{\cal
E})$-valued function in $\zeta = \rho + i \delta$,
given by
\[  \big(V(\varphi) a\big)(\rho + i \delta, \eta) =
       \big(V(\varphi_\delta) a\big)(\rho, \eta)     \]
for $\varphi_\delta(\theta) := e^{\theta \delta}
\varphi(\theta)$. From Theorem \ref{734.cutR0703.th} we
obtain $\big(V(\varphi) a\big)(\rho + i \delta, \eta) \in
S^\mu(\R^{1+q}; {\cal E}, \widetilde{\cal E})$ for every
$\delta \in \R$. By virtue of the continuity of $\delta
\to \varphi_\delta$, $\R \to C_0^\infty(\R)$ and of the
continuity of \eqref{734.cont0703.eq} it follows that
\eqref{734.conthol803.eq} induces a continuous mapping
\[  V : C_0^\infty(\R) \to S^\mu(\R^{1+q}; {\cal E},
        \widetilde{\cal E}) \to S^\mu (I_\delta \times
	\R^q; {\cal E},\widetilde{\cal E}),    \]
$I_\delta := \{ \zeta \in \C : \imb \zeta = \delta \}$
wich is uniform in compact $\delta$-intervals. The
closed graph theorem gives us also the continuity of
\eqref{734.conthol803.eq} with respect to the Fr\'echet
topology of $S^\mu(\boldsymbol{C} \times \R^q; {\cal E},
\widetilde{\cal E})$.	      
\qed \\

\subsection{Meromorphic Mellin symbols and operators with asymptotics}

As an ingredient of our cone algebra we now study meromorphic Mellin symbols, starting from $M_\mathcal{O}^\mu(\E,\widetilde\E)$ (see Definition \ref{defmo} for $q=0$).
\begin{Thm} \label{thm331}
$h\in M_\mathcal{O}^\mu (\E,\widetilde\E)$ and $h|_{\Gamma_\beta}\in S^{\mu-\varepsilon}(\Gamma_\beta;
\E,\widetilde\E)$ for some $\varepsilon>0$ entails $h\in M_\mathcal{O}^{\mu-\varepsilon}(\E,
\widetilde\E)$.
\end{Thm}
\proof
The ideas of the proof are similar to the case of the cone calculus with smooth base $X$ and the scales $ \big(H^s(X)\big)_{s\in\R}$ (see, e.g., the thesis of Seiler \cite{Seil3}). 
\qed
\begin{Prop} \label{mmult}
Let $h(w)\in M_\mathcal{O}^\mu(\E_0,\widetilde\E)$, $f(w)\in M_\mathcal{O}^\nu(\E,\E_0)$; then for pointwise composition we have $h(w)f(w)\in M_\mathcal{O}^{\mu+\nu}(\E,\widetilde\E)$.
\end{Prop}
\proof
The proof is obvious.
\qed
\begin{Def} \label{Def7347}
An element $h(w)\in M^\mu_\mathcal{O}(\E,\widetilde\E)$ is called elliptic, if for some $\beta\in
\R$ the operators $h(\beta+i\rho):E^s\to\widetilde E^{s-\mu}$ are invertible for all $s\in\R$ and
$h^{-1}(\beta+i\rho)\in S^{-\mu}(\R_\rho;\widetilde\E,\E)$.
\end{Def}
\begin{Thm} \label{Thm7348}
Let $h\in M^\mu_\mathcal{O}(\E,\widetilde\E)$ be elliptic. Then,
\begin{equation}\label{opesmu}
h(w):E^s\to\widetilde E^{s-\mu}
\end{equation}
is a holomorphic family of Fredholm operators of index zero for $s\in\R$. There is a set 
$D\subset\C$, with $D\cap\{c\leq\reb\ w\leq c'\}$ finite for every $c\leq c'$, such that the operators
\eqref{opesmu} are invertible for all $w\in\C\setminus D$.
\end{Thm}
\proof
By assumption we have $g:=(h|_{\Gamma_\beta})^{-1}\in S^{-\mu}(\Gamma_\beta;\widetilde\E,\E)$. 
Applying a version of the kernel cut-off construction, now referring to parallels of the imaginary axis rather than the real axis, with a function $\psi\in C_0^\infty(\R_+)$, $\psi\equiv 1$ near $1$, we 
obtain a continuous operator
\[ V(\psi ):S^{-\mu }(\Gamma_\beta ;\widetilde\E ,\E )\to M^{-\mu }_\mathcal{O}(\widetilde\E ,\E )\]
where $V(\psi)g|_{\Gamma_\beta}=g\ \mod\ S^{-\infty}(\Gamma_\beta;\widetilde\E,\E)$. Setting 
$h^{(-1)}(w):=V(\psi)g$ we obtain $h^{(-1)}(w)\in M^{-\mu}_\mathcal{O}(\widetilde\E,\E)$, and from 
proposition \ref{mmult} it follows that
\[ h(w)h^{(-1)}(w)\in M^0_\mathcal{O}(\widetilde\E,\widetilde\E) , \
h(w)^{(-1)}h(w)\in M^0_\mathcal{O}(\E ,\E ) \]
and 
\begin{equation} \label{relathh-1}
h(w)h^{(-1)}(w)|_{\Gamma_\beta}-1\in S^{-\infty}(\Gamma_\beta ;\widetilde\E ,\widetilde\E ), 
h(w)^{(-1)}h(w)|_{\Gamma_\beta}-1\in S^{-\infty}(\Gamma_\beta ;\E ,\E ). 
\end{equation}
for every $\beta\in\R$, and hence
\begin{equation}\label{hencehh-1}
h(w)h^{(-1)}(w)=1+m(w),\ h(w)^{(-1)}h(w)=1+l(w)
\end{equation} 
for certain $m(w)\in M^{-\infty}_\mathcal{O}(\widetilde\E,\widetilde\E)$, $l(w)\in M^{-\infty}
_\mathcal{O}(\E,\E)$. For every $s\in\R$ and every fixed $w\in\C$ the operators
\[ m(w):\widetilde E^s\to \widetilde E^\infty ,\ l(w):E^s\to E^\infty \]
are continuous. Therefore, since the scales have the compact embedding property, from 
\eqref{hencehh-1} we obtain that $h^{(-1)}(w)$ is a two-sided parametrix of $h(w)$ for every $w$,
i.e., the operators \eqref{opesmu} are Fredholm. Since $h(w)\in\mathcal{A}(\C,\Li^\mu(E^s,
E^{s-\mu}))$ is continuous in $w\in\C$ we have $\ind\ h(w_1)=\ind\ h(w_2)$ for every $w_1,w_2\in\C$.
However, since $h(w)=0$ consists of invertible operators on the line $\Gamma_\beta$ it follows that 
$\ind\ h(w)=0$ for all $w\in\C$. Finally, from the realtions \eqref{opesmu} we see that for every
$c\leq c'$ there is an $L(c,c')>0$ such that the operators \eqref{opesmu} are invertible for all
$w\in\C$ with $|\imb\ w|\geq L(c,c'),\ c\leq\reb\ w\leq c'$. Then a general result on holomorphic 
Fredholm families tells us that the strip $c\leq\reb\ w\leq c'$ contains at most finitely many 
points where \eqref{opesmu} is not invertible. Those points just constitute the set $D$, it is also 
independent of $s$, since $\textup{ker}\ h(w)$ is independent of $s$ as we easily see from 
\eqref{hencehh-1} and the smoothing remainders; then vanishing of the index shows that the 
invertibility holds exactly when $\textup{ker}\ h(w)=0$.
\qed
\begin{Thm}
The ellipticity of $h$ with respect to $\Gamma_\beta$ as in Definition \textup{\ref{Def7347}} entails the ellipticity with respect to $\Gamma_\delta$ for all $\delta\in\R$ such that $\Gamma_\delta \cap D =\emptyset$. In that sense Definition \textup{\ref{Def7347}} is independent of the choice of $\beta$.
\end{Thm}
\proof
Let us apply the kernel cut-off operator $V(\psi_\varepsilon)$, where $\psi_\varepsilon\in C_0^\infty (\R_+)$ is of the form $\psi_\varepsilon(t)=\psi(\varepsilon t)$, $\varepsilon>0$, for some cut off fuction $\psi$. Then, setting
\[ V(\psi_\varepsilon)(h^{-1}(\beta+i\rho))=:f_\varepsilon\in M^{-\mu}_\mathcal{O}(\widetilde\E,\E) \]
we have $f_\varepsilon|_{\Gamma_\beta}\in S^{-\mu}(\Gamma_\beta;\widetilde\E,\E)$ and $f_\varepsilon| _{\Gamma_\beta}\to h^{-1}(\beta+i\rho)$ as $\varepsilon\to 0$ in the topology of $S^{-\mu}(\Gamma_\beta;\widetilde\E,\E)$. This shows us that $f_{\varepsilon_1}|_{\Gamma_\beta}$ is pointwise invertible when $\varepsilon_1>0$ is sufficiently small. Let us set $h^{(-1)}(w)= f_{\varepsilon_1}(w)$. According to Proposition \ref{mmult} we have $g(w):=h^{(-1)}(w)h(w)\in M_\mathcal{O}^0(\E,\E)$ and by construction 
\[ g|_{\Gamma_\beta}=1+l\ \textup{for some}\ l\in S^{-\infty}(\Gamma_\beta;\E,\E). \]
Then Theorem \ref{thm331} yields $g=1\mod M_\mathcal{O} ^{-\infty}(\E,\E)$. It follows that
\[ h^{-1}|_{\Gamma_\delta}h|_{\Gamma_\delta}=1+l_\delta\ \textup{for some}\ l_\delta\in S^{-\infty} (\Gamma_\delta;\E,\E) \]
and hence, since $h|_{\Gamma_\delta}$ is pointwise invertible, 
\[ h^{(-1)}|_{\Gamma_\delta}=(1+l_\delta)\left(h|_{\Gamma_\delta}\right)^{-1}. \]
which yields
\begin{equation}
\left(h|_{\Gamma_\delta}\right)^{-1}=(1+l_\delta)^{-1}h^{(-1)}|_{\Gamma_\delta}.
\end{equation}
From Proposition \ref{fs-} we know that $l_\delta\in\Sch(\Gamma_\delta,\Li^{-\infty}(\E,\E))$ and it is also clear that $(1+l_\delta)^{-1}=1+m_\delta$ for some $m_\delta\in  \Sch(\Gamma_\delta, \Li^{-\infty}(\E,\E))$. Then Proposition \ref{s,s-nu} shows that $(h|_{\Gamma_\delta})^{-1}\in S^{-\mu}(\Gamma_\delta;\widetilde\E,\E)$.
\qed \\
\indent A sequence
\[ R=\{(p_j,m_j,L_j)\}_{j\in\Z} \]
is called a discrete asymptotic type of Mellin symbols, if $p_j\in\C$, $m_j\in\N$, and $L_j\subset \Li^{-\infty}(\E,\widetilde\E)$ is a finite-dimensional subspace of finite  rank operators; moreover, $\pi_\C R:=(p_j)_{j\in\Z}$ is assumed to intersect the strips $\{w\in\C:c_1\leq\reb w\leq c_2\}$ in a finite set, for every $c_1\leq c_2$. Let $M^{-\infty}_R(\E,\widetilde\E)$ denote the space of all functions $m\in\An(\C\setminus\pi_\C R,\Li^{-\infty}(\E,\widetilde\E))$ which are meromorphic with poles at the points $p_j$ of multiplicity $m_j+1$ and Laurent coefficients at $(w-p_j)^{-(k+1)}$ belonging to $L_j$ for $0\leq k\leq m_j$, and $\chi(w)m(w)|_{\Gamma_\delta}\in \Sch(\Gamma_\delta; \E,\widetilde\E)$ for every $\delta\in\R$, uniformly in compact $\delta$-intervals, where $\chi$ is any $\pi_\C R$-excision function. Moreover, we set
\begin{equation} \label{neu34} M^\mu_R(\E,\widetilde\E):=M^\mu_\mathcal{O}(\E,\widetilde\E)+M^{-\infty}_R(\E,\widetilde\E) .
\end{equation}

\begin{Thm} \label{mulom}
Let $h\in M^\mu_R(\E_0,\widetilde\E)$, $f\in M^\nu_S(\E,\E_0)$ with asymptotic types $R,S$ and orders $\mu,\nu\in\R$, then we have $hf\in M^{\mu+\nu}_P(\E,\widetilde\E)$ with some resulting asymptotic type $P$.
\end{Thm}
\proof
The proof of this result is analogous to the one in the ``concrete'' cone calculus, see \cite{Schu2}.
\qed
\begin{Prop} \label{invo}
For every $m\in M^{-\infty}_R(\E,\E)$ there exists an $m^{(-1)}\in M^{-\infty}_S(\E,\E)$  with another asymptotic type $S$ such that 
\[ \big(1+m(w)\big)\big(1+m^{(-1)}(w)\big)=1 . \]
\end{Prop}
For the proof we employ the following Lemma.
\begin{Lem} \label{lem378}
Let $E$ be a Banach space, $U\subseteq\C$ open, $0\in U$, and let $h\in\An(U,\Li(E))$ be an element such that $h(w)=0$ for all $u\in F$ where $F\subseteq E$ is a closed subspace of finite codimension. Moreover, let $a_1,\dots,a_N\in\Li(E)$ be operators of finite rank, for some $N\in\N\setminus\{0\}$. Then there is a $\delta>0$ such that the meromorphic $\Li(E)$-valued function
\[ f(w)=1+h(w)+\sum_{j=1}^N a_j w^{-j} \]
is invertible for all $0<|w|<\delta$.
\end{Lem}
\proof[\textup{\textbf{Proof of Proposition \ref{invo}}}]
First observe that if $m\in\Li^{-\infty}(\E,\E)$ is an operator such that 
\[ 1+m:E^s\to E^s \]
is invertible for all $s\in\R$, we can define an operator $g\in\Li^0(\E,\E)$ such that $(1+m)(1+g)=1$. This gives us $1+m+g+mg=1$, and $m,mg\in\Li^{-\infty}(\E,\E)$ implies $g\in\Li^{-\infty}(\E,\E)$. \\
Moreover, our operator function $1+m$ is holomorphic in $\C\setminus\pi_\C R$. Then $g=(1+m)^{-1}-1$ is holomorphic in $\C\setminus D$ with values in $L^{-\infty}(\E,\E)$, where $D\subseteq\C$ is a countable set such that $\{w\in\C:c_1\leq \reb w\leq c_2\}\cap \{w\in\C:\textup{dist}(w,\pi_\C R)>\varepsilon \}\cap D$ is finite for every $c_1\leq c_2$ and $\varepsilon>0$. If $\chi(w)$ is any $D$-excision function, then we have
\[\chi(w)(1+m(w))|_{\Gamma_\beta}\in\Sch(\Gamma_\beta,\Li^{-\infty}(\E,\E)) \]
for every $\beta\in\R$, uniformly in compact $\beta$-intervals. This shows that $\chi(w)(1+m(w))^{-1} |_{\Gamma_\beta}\in\Sch(\Gamma_\beta,\Li^{-\infty}(\E,\E))$ for every $\beta\in\R$, uniformly in compact $\beta$-intervals. It remains to show that $g$ is meromorphic with poles at the points of $D$, that $D$ has no accumulation points at $\pi_\C R$, and that the Laurent coefficients are of the desired kind, namely, to belong to $L^{-\infty}(\E,\E)$ and to be of finite rank. Let us verify that there are no accumulation points of the singularities of $(1+m(w))^{-1}$. let $w_0$ be a pole of $m$, i.e., $w_0\in\pi_\C R$. Then we can write
\[ 1+m(w)=1+m_0(w)+\sum_{k=1}^K b_k(w-w_0)^{-k} \]
with suitable $K\in\N$, $m_0$ holomorphic in a neighbourhood of $w_0$ and $\Li^{-\infty} (\E,\E)$-valued, with finite rank operators $b_k$. Note that $m_0\not\equiv -1$. Setting $n(w):= \sum_{k=1}^K b_k(w-w_0)^{-k}$ we have
\[ 1+m(w)=(1+m_0(w))(1+(1+m_0(w))^{-1}n(w)). \]
Since $m_0$ is holomorphic near $w_0$ and $1+m_0(w)$ a Fredholm family, the singularities of $(1+m_0(w))^{-1}$ form a countable discrete set; therefore there is a $\delta>0$ such that $(1+m_0(w))^{-1}$ exists for all $0<|w-w_0|<\delta$. Moreover, $(1+m_0(w))^{-1}n(w)$ can be written in the form $h(w)+\sum_{j=1}^N a_j(w-w_0)^{-j}$ with a suitable $h$ which is holomorphic near $w_0$ and finite rank operators $a_j$, $1\leq j\leq N$. The operator $(1+m_0(w))^{-1}n(w)$ vanishes on the space $F:=\bigcap_{k=1}^K \textup{ker}\, b_k$ which is of finite codimension. Setting $M:=\bigcap_{j=1}^N \textup{ker}\, a_j$ it follows that $h(w)u=0$ for all $u\in M\cap F$; the latter space is also of finite codimension. Lemma \ref{lem378} then shows that $1+(1+m_0(w))^{-1}n(w)$ is invertible in $0<|w-w_0|<\delta$ for a suitable $\delta>0$.
\qed 
\begin{Thm}
Let $h\in M^\mu_\mathcal{O}(\E,\widetilde\E)$ be elliptic, then there is an $f\in M^{-\mu}_S (\widetilde\E,\E)$ with asymptotic type $S$ such that $hf=1$.
\end{Thm}
\proof
Let $h^{(-1)}(w)\in M_\mathcal{O}^\mu(\E,\widetilde\E)$ be as in the proof of Theorem \ref{Thm7348}. Then we have the relation \eqref{hencehh-1}. By virtue of Proposition \ref{invo} there exists a $g\in M_P^\infty(\E,\E)$ for some asymptotic type $P$ such that $(1+m(w))(1+g(w))=1$. This yields $h(w)f(w)=1$ for $f:=h^{(-1)}(1+g)$ which belongs to $M_S^{-\mu}(\widetilde\E,\E)$, according to Theorem \ref{mulom}. In a similar manner we find an $\td f\in M_{\widetilde S}^{-\mu}(\widetilde\E, \E)$ such that $\td f(w)h(w)=1$. This implies $f=\td f$.
\qed 
\begin{Def}
A $g\in M^\mu_R (\E,\widetilde\E)$ is said to be elliptic, if there is a $\beta\in\R$ such that $(g|_{\Gamma_\beta})^{-1}\in S^{-\mu}(\R_\rho;\widetilde\E,\E)$.
\end{Def}
\begin{Thm}
If $g\in M^\mu_R (\E,\widetilde\E)$ is elliptic, there is an $f\in M^{-\mu}_S(\widetilde\E,\E)$ such that $gf=1$.
\end{Thm}
\proof
Applying the kernel cut-off operator to $(g|_{\Gamma_\beta})^{-1}$ we find an $h^{(-1)}\in M^{-1}_\mathcal{O}(\widetilde\E,\E)$ such that $h^{(-1)}|_{\Gamma_\beta}-(g|_{\Gamma_\beta})^{-1}\in S^{-\infty}(\Gamma_\beta;\widetilde\E,\E)$. By definition we have $g=g_0+g_1$ for certain $g_0\in M^\mu_\mathcal{O}(\E,\widetilde\E)$, $g_1\in M^{-\infty}_R(\E,\widetilde\E)$. Then $h^{(-1)}g_0|_{\Gamma_\beta}=1 \mod S^{-\infty}(\Gamma_\beta;\E,\E)$ implies $h^{(-1)}g_0=1 \mod M_\mathcal{O}^\infty(\E,\E)$ (see Theorem \ref{thm331}). If follows that $h^{(-1)}g=1+m$ for some $m\in M_R^{-\infty}(\E,\E)$ with an aymptotic type $R$ (see Theorem \ref{mulom}). Thus Proposition \ref{invo} gives us $g^{-1}=(1+m)^{-1}h^{(-1)}\in M_S^{-\mu}(\widetilde\E,\E)$ with some asymptotic type $S$.
\qed \\
\indent Parallel to the spaces of Mellin symbols \eqref{neu34} we now introduce subspaces of $\Hi^{s,\gamma}(\R_+,\E)$ with discrete asymptotics. We consider a sequence 
\begin{equation} 
\label{Pasympeq}
P:=\{(p_j,m_j)\}_{0\leq j\leq N}
\end{equation}
with $N\in\N\cup\{+\infty\}$, $m_j\in\N$, $0\leq j\leq N$.\\
A sequence \eqref{Pasympeq} is said to be a discrete asymptotic type, associated with weight data $(\gamma,\Theta)$ (with a weight $\gamma\in\R$ and a weight interval $\Theta=(\vartheta,0]$, $-\infty\leq\vartheta\leq 0$), if 
\[ \pi_\C P:=\{p_j\}_{0\leq j\leq N}\subset \{ w\in\C:\frac{d+1}{2}-\gamma+\vartheta<\reb\, w< \frac{d+1}{2}-\gamma\}, \]
and $\pi_\C P$ is finite when $\vartheta$ is finite, and $\reb\, p_j\to -\infty$ as $j\to\infty$ when $\vartheta =-\infty$ and $N=+\infty$. We will say that $P$ satisfies the shadow condition, if $(p,m)\in P$ implies $(p-j,m)\in P$ for all $j\in\N$ with $\frac{d+1}{2}-\gamma+\vartheta< \reb\, (p-j)<\frac{d+1}{2}-\gamma$.\\
If $\Theta$ is finite we define the (finite-dimensional) space
\[ S_P(\R_+,\E):=\Big\{\sum_{j=0}^N\sum_{k=0}^{m_j}\omega(r)c_{jk}r^{-p_j}\log^k r:\ c_{jk}\in E^\infty,\  0\leq k\leq m_j,\ 0\leq j\leq N\Big\} \]
with some fixed cut-off function $\omega$ on the half-axis. We then have
\[S_P(\R_+,\E)\subset\Hi^{\infty,\gamma}(\R_+,\E). \]
Moreover, we set 
\[ \Hi_\Theta^{s,\gamma}(\R_+,\E):=\lim_{\overleftarrow{m\in\N}}\{ \omega\Hi^{s,\gamma-\vartheta -\frac{1}{m+1}}(\R_+,\E)+(1-\omega)\Hi^{s,\gamma}(\R_+,\E)\} \]
endowed with the Fr\'echet topology of the projective limit, and 
\[ \Hi_P^{s,\gamma}(\R_+,\E):=\Hi^{s,\gamma}_\Theta(\R_+,\E) + S_P(\R_+,\E) \]
as a direct sum of Fr\'echet spaces.\\
In order to formulate spaces with discrete asymptotics of type $P$ in the case $\Theta=(-\infty,0]$ we from $P_k:=\{(p,m)\in P:\reb\, p>\frac{d+1}{2}-\gamma-(k+1)\}$ for any $k\in\N$. From the above constuction we have the spaces $\Hi_{P_k}^{s,\gamma}(\R_+,\E)$ together with continuous embeddings 
\[ \Hi^{s,\gamma}_{P_{k+1}}(\R_+,\E)\hookrightarrow \Hi_{P_k}^{s,\gamma}(\R_+,\E),\ k\in\N. \]
We then define
\begin{equation} \label{projlimas}
\Hi_P^{s,\gamma}(\R_+,\E):=\lim_{\overleftarrow{k\in\N}}\Hi_{P_k}^{s,\gamma}(\R_+,\E)
\end{equation}
in the corresponding Fr\'echet topology of the projective limit.
\begin{Rem}
The relation $u\in\Hi_P^{s,\gamma}(\R_+,\E)$ with $P$ being associated with $(\gamma,\Theta)$, $\Theta=(-\infty,0]$, is equivalent with the existence of \lr unique\rr\ coefficients $c_{jk}\in E^\infty$, $0\leq k\leq m_j$, such that for every $l\in\R_+$ there is an $N=N(l)\in\N$ with
\[ \omega(r)\Big(u(r,x)-\sum_{j=0}^N\sum_{k=0}^{m_j}c_{jk}r^{-p_j}\log^k r\Big)\in\Hi^{s,\gamma+l} (\R_+,\E). \]
\end{Rem}
Similarly as in the ``concrete'' cone calculus (see \cite{Schu2}) we have the following continuity result:
\begin{Thm}
For every $f\in M^\mu_R(\E,\widetilde\E)$ the operator \eqref{733.cof2201.eq} restricts to a continuous operator
\[ \sop_M^\gamma(f):\Hi_P^{s,\gamma}(\R_+,\E)\to\Hi_Q^{s-\mu,\gamma}(\R_+,\widetilde\E) \]
for every $s\in\R$ and every asymptotic type $P$ with some resulting $Q$.
\end{Thm}
The case  of Mellin symbols with variable coefficients is also of interest in the corner calculus. It is then adequate to assume $f(r,w)\in C^\infty(\Rb_+,M^\mu_R(\E,\widetilde\E))$ and to consider operators $\omega \sop_M^\gamma(f)\td\omega$ in combination with cut-off functions $\omega(r),\td\omega(r)$. Those induce continuous operators $\Hi_P^{s,\gamma}(\R_+,\E)\to\Hi_Q^{s-\mu,\gamma}(\R_+,\widetilde\E)$ as well.

\end{document}